\documentclass[11pt]{article}

\usepackage[margin=1in]{geometry}
\usepackage{amsmath,amssymb,amsthm,mathtools}
\usepackage{microtype}
\usepackage{needspace,etoolbox}
\usepackage[colorlinks=true,linkcolor=blue,
citecolor=red,urlcolor=blue]{hyperref}
\usepackage[nameinlink,capitalize,noabbrev]{cleveref}

\newtheorem{theorem}{Theorem}[section]
\newtheorem{lemma}[theorem]{Lemma}
\newtheorem{proposition}[theorem]{Proposition}
\newtheorem{corollary}[theorem]{Corollary}
\numberwithin{equation}{section}

\newcommand{\R}{\mathbb R}
\newcommand{\dd}{\,\mathrm d}

\newcommand{\C}{\mathbb C}
\newcommand{\Hh}{\mathcal H}
\newcommand{\Dcal}{\mathcal D}
\newcommand{\Ncal}{\mathcal N}
\newcommand{\Fcal}{\mathcal F}
\newcommand{\Lcal}{\mathcal L}
\newcommand{\ip}[2]{\langle #1,#2\rangle}
\DeclareMathOperator{\dist}{dist}
\DeclareMathOperator{\spanop}{span}
\DeclareMathOperator{\ind}{ind}
\DeclareMathOperator{\Rea}{Re}
\allowdisplaybreaks[2]
\BeforeBeginEnvironment{theorem}{\Needspace{5\baselineskip}}
\BeforeBeginEnvironment{lemma}{\Needspace{5\baselineskip}}
\BeforeBeginEnvironment{proposition}{\Needspace{6\baselineskip}}
\BeforeBeginEnvironment{corollary}{\Needspace{5\baselineskip}}

\title{Ground state uniqueness and quantitative stability for the mass-critical Hartree equation}

\author{Guangze Gu and Zhipeng Yang\thanks{Corresponding author: yangzhipeng326@163.com}}

\date{}

\AtEndDocument{
\par
\bigskip
\bigskip

\noindent
\textbf{Guangze Gu}\\[0.2em]
\textsc{Department of Mathematics, Yunnan Normal University, Kunming 650500, China}\\
\textsc{Yunnan Key Laboratory of Modern Analytical Mathematics and Applications, Kunming 650500, China}\\[0.3em]
\textit{E-mail address}: \texttt{guangzegu@163.com}

\bigskip

\noindent
\textbf{Zhipeng Yang}\\[0.2em]
\textsc{Department of Mathematics, Yunnan Normal University, Kunming 650500, China}\\
\textsc{Yunnan Key Laboratory of Modern Analytical Mathematics and Applications, Kunming 650500, China}\\[0.3em]
\textit{E-mail address}: \texttt{yangzhipeng326@163.com}
}

\begin{document}
\maketitle

\begin{abstract}
We prove that the ground state of the mass-critical Hartree equation is unique up to phase rotations and translations and nondegenerate in every dimension $d\ge3$. This makes the recent conditional $L^2$ classification of minimal-mass blow-up and non-scattering solutions unconditional; we also give a direct $H^1$ proof of the finite-time blow-up classification. We further establish a global quadratic stability estimate for the sharp Hartree interpolation inequality. For the optimal stability constant, we rule out vanishing of minimizing sequences and describe their possible loss of compactness.
\end{abstract}



\Needspace{7\baselineskip}
\section{Introduction and main results}\label{Sec1}

We consider the stationary mass-critical Hartree equation
\begin{equation}\label{eq1.1}
-\Delta Q+Q=(|x|^{-2}*|Q|^2)Q
\qquad\text{in }\R^d,\qquad d\ge3,
\end{equation}
whose solutions give the solitary waves $e^{it}Q$ of the focusing Hartree flow
\begin{equation}\label{eq1.2}
i\partial_tu+\Delta u+(|x|^{-2}*|u|^2)u=0,
\qquad u(0)=u_0.
\end{equation}
The restriction $d\ge3$ is natural because $|x|^{-2}$ is not locally integrable in dimensions one and two. Ground states of \eqref{eq1.1} optimize a sharp nonlocal interpolation inequality, and their common mass is the threshold for global existence and scattering of \eqref{eq1.2} \cite{MXZ2009,MMRZ2026}. For the local mass-critical Schr\"odinger equation, Weinstein's variational characterization \cite{Weinstein1983} and Kwong's uniqueness theorem \cite{Kwong1989} are the elliptic basis of the threshold theory, beginning with Merle's classification of minimal-mass blow-up \cite{Merle1993}. For \eqref{eq1.1}, Krieger, Lenzmann and Rapha\"el \cite{KLR2009} proved uniqueness and radial nondegeneracy in dimension four, where $|x|^{-2}$ is the Newtonian kernel, by adapting the argument of Lieb \cite{Lieb1977}. For $d\ne4$, uniqueness has been assumed in several results on threshold dynamics \cite{LiZhang2009,GaoWang2017,MMRZ2026,MMZ2026}; Li and Zhang explicitly record the unresolved higher-dimensional case in \cite{LiZhang2009}. In this paper we prove that the ground state is unique and nondegenerate in every dimension $d\ge3$, which removes this hypothesis from those results.

After radial reduction, \eqref{eq1.1} is still an integro-differential equation, so the shooting and comparison arguments behind Kwong's theorem are not directly available. Newton's theorem turns the potential into the solution of a local Poisson equation only when $d=4$:
\[
\int_{\mathbb S^3}\frac{\dd\theta}{|re_1-s\theta|^2}
=|\mathbb S^3|\max\{r,s\}^{-2},
\qquad
\int_{\mathbb S^2}\frac{\dd\theta}{|re_1-s\theta|^2}
=\frac{2\pi}{rs}\log\frac{r+s}{|r-s|}.
\]
For $d\ne4$ the spherical average depends nontrivially on both radii, and for $d=3$ it is logarithmically singular on the diagonal. Denote by $a_\nu(r,s)$ this average with the dimension replaced by a real parameter $\nu\ge3$, as defined in \eqref{eq2.1}. Newton's identity gives a vanishing mixed derivative off the diagonal at $\nu=4$. In every real dimension, Lemma~\ref{Lem2.1} gives
\[
B_\nu(r,s)=\partial_r\partial_sa_\nu(r,s)
=2(\nu-4)k_\nu(r,s),\qquad k_\nu(r,s)>0,\qquad r\ne s.
\]
Its sign drives two comparison arguments, one above and one below dimension four, which allow us to continue the four-dimensional ground state in $\nu$.

Throughout the paper, the convolution kernel is the unnormalized function $|x|^{-2}$. For $u\in H^1(\R^d;\C)$, write
\[
M(u)=\int_{\R^d}|u|^2,\qquad
A(u)=\int_{\R^d}|\nabla u|^2,\qquad
\Dcal(u)=\iint_{\R^d\times\R^d}
\frac{|u(x)|^2|u(y)|^2}{|x-y|^2}\dd x\dd y,
\]
and set
\[
E(u)=\frac12A(u)-\frac14\Dcal(u),\qquad
S(u)=E(u)+\frac12M(u).
\]
These functionals are finite on $H^1$ by \eqref{eq1.3} below. A ground state is a solution $Q\in H^1(\R^d;\C)\setminus\{0\}$ of \eqref{eq1.1} with least action $S$ among all such solutions.

At a positive ground state of \eqref{eq1.1}, the linearized operators for real and imaginary perturbations are
\[
\begin{aligned}
L_+h&=-\Delta h+h-V_Qh
-2Q\bigl(|x|^{-2}*(Qh)\bigr),\\
L_-h&=-\Delta h+h-V_Qh,
\qquad V_Q=|x|^{-2}*Q^2.
\end{aligned}
\]
Their Morse index is the number of negative eigenvalues, counted with multiplicity, and their kernels are understood in the weak $H^1(\R^d;\R)$ sense. We write $q_\pm[h]=\ip{L_\pm h}{h}$ for their quadratic forms. We use $\ip{f}{g}$ for the real $L^2$ pairing on the domain in question, extended to the corresponding duality pairing when one argument lies in the dual of the energy space, and
\[
\ip{u}{v}_{H^1}=\int_{\R^d}
\bigl(u\overline v+\nabla u\cdot\nabla\overline v\bigr)
\]
for the complex $H^1$ inner product.

\begin{samepage}
\begin{theorem}\label{Thm1.1}
For every integer $d\ge3$, equation \eqref{eq1.1} has a unique positive radial ground state $Q_d$. Every complex ground state is of the form
\[
e^{i\theta}Q_d(\cdot-y),\qquad \theta\in\R,\quad y\in\R^d.
\]
Moreover,
\[
\ker L_+=\spanop\{\partial_{x_1}Q_d,\ldots,\partial_{x_d}Q_d\},
\qquad
\ker L_-=\spanop\{Q_d\}.
\]
The real operator $L_+$ has Morse index one. If $d\ge4$, every positive real $H^1$ solution of \eqref{eq1.1} is a translate of $Q_d$.
\end{theorem}
\end{samepage}

For $d=4$, uniqueness and radial nondegeneracy are due to Krieger, Lenzmann and Rapha\"el \cite[Theorems~3--4]{KLR2009}, and the full kernel characterization is given in \cite[Lemma~2.5]{GSW2026}. The new cases in Theorem~\ref{Thm1.1} are $d\ne4$. For $d\ge4$, the theorem includes all positive $H^1$ solutions. Indeed, Proposition~\ref{Prop3.3} makes radial Morse index one automatic for positive decreasing profiles above dimension four, and the symmetry theorem of Ma and Zhao \cite{MaZhao2010} reduces positive solutions to this class. For $d=3$, Proposition~\ref{Prop4.3} gives uniqueness up to translation among positive solutions with radial Morse index one, without an action-minimizing assumption. The uniqueness, up to translation, of positive $H^1$ solutions in dimension three remains an open problem.

To prove Theorem~\ref{Thm1.1}, we embed the radial equation into a family indexed by the real dimension $\nu\ge3$. The conjugation $z(r)=r^{(\nu-1)/2}Q(r)$ places this family on $H_0^1(0,\infty)$. We then combine radial nondegeneracy and compactness to continue uniqueness from $\nu=4$, following the continuation strategy of Frank and Lenzmann \cite{FL2013} and Frank, Lenzmann and Silvestre \cite{FLS2016}.

The mixed derivative $B_\nu$ enters the conjugated radial quadratic form through a cumulative change of variable. For $y\perp z$, set
\[
F_y(r)=\int_0^r z(s)y(s)\dd s,\qquad
y=F_y'/z,\qquad \rho=z^2.
\]
In this variable, Lemma~\ref{Lem3.2} writes the local part of $q_\nu$ as $\int_0^\infty\rho\,|(\rho^{-1}F_y')'|^2\dd r$ and the nonlocal part as $-2\ip{B_\nu}{F_y\otimes F_y}$. For $\nu>4$, $B_\nu$ is a positive locally integrable kernel. The cumulative function $w=F_Z$ of the dilation mode is positive. Inverting the local operator twice through its positive Green kernel gives the comparison with $w$ in \eqref{eq3.2}; its equality case leaves only multiples of $Z$. For $3\le\nu<4$, $B_\nu$ is a distribution with a negative off-diagonal part. Pairing it with a squared difference removes the diagonal singularity and gives \eqref{eq3.6}, whose nonlocal term is $-\iint\mathcal K_z(t,u)k(t)k(u)\dd t\dd u$ with $\mathcal K_z>0$ and $k=(F_y/w)'$. Replacing $k$ by $|k|$ leaves the local terms unchanged and strictly decreases the quadratic form whenever $k$ changes sign. Radial Morse index one makes the form nonnegative on $z^\perp$. Hence a kernel element has monotone $F_y/w$, and $\ip{Vz}{y}=0$ together with $V'<0$ forces it to vanish. At $\nu=4$, we use the uniqueness and radial nondegeneracy theorems of \cite{KLR2009}.

The comparison function $w=F_Z$ comes from mass criticality. For the dilation derivative $Z=z/2+rz'$ and the radial linearization $\Lcal_\nu$,
\[
\Lcal_\nu Z=-2z,\qquad \ip{z}{Z}=0,\qquad
F_Z(r)=\frac12rz(r)^2>0.
\]
The last function vanishes at both endpoints. Together with radial Morse index one, the first two identities give nonnegativity on the mass tangent space. These facts provide the comparison function above dimension four and the cancellation below it. The continuation method can be used for other interactions once compactness and invertibility are available, but the comparison also needs the inverse-square kernel identity and a frequency dilation tangent to the mass constraint. For other homogeneities, the cumulative function need not vanish at infinity.

Proposition~\ref{Prop2.4} gives uniform bounds on compact parameter intervals for all positive decreasing profiles, without assuming that they minimize the action. An annular mass estimate controls the radius at which the potential reaches $1/4$, and an exterior comparison bounds the total mass. Radial nondegeneracy makes the number of profiles with radial Morse index one finite and locally constant. Uniqueness at $\nu=4$ fixes this number at one throughout the parameter range.

To identify the full Euclidean kernel, we decompose in spherical harmonics. The angular interaction kernels strictly decrease with the harmonic degree by Lemma~\ref{Lem4.5}. Translation modes account for degree one, and the higher degrees are strictly positive.

The Hardy--Littlewood--Sobolev inequality \cite{Lieb1983} and Sobolev interpolation give
\begin{equation}\label{eq1.3}
\Dcal(u)\le C_{d,1}\|u\|_{2d/(d-1)}^4
\le C_{d,2}M(u)A(u).
\end{equation}
The corresponding interpolation problem is
\[
j_d=\inf_{u\in H^1(\R^d;\C)\setminus\{0\}}
\frac{M(u)A(u)}{\Dcal(u)},\qquad m_d=2j_d.
\]
The infimum is positive and attained \cite{MorozVS2013}. The variational characterization and the Pohozaev identity give, for every ground state,
\begin{equation}\label{eq1.4}
M(Q)=A(Q)=m_d,\qquad \Dcal(Q)=2m_d,\qquad E(Q)=0,
\end{equation}
as recalled in Lemma~\ref{Lem4.1}. The sharp form of \eqref{eq1.3} is
\begin{equation}\label{eq1.5}
\Dcal(u)\le\frac{2}{m_d}M(u)A(u).
\end{equation}

Theorem~\ref{Thm1.1} removes the uniqueness and nondegeneracy assumptions from several results on the threshold dynamics of \eqref{eq1.2}. This flow is invariant under $u(t,x)\mapsto\lambda^{d/2}u(\lambda^2t,\lambda x)$ and conserves mass and energy along $H^1$ solutions. By \eqref{eq1.5},
\[
E(u)\ge\frac12\left(1-\frac{M(u)}{m_d}\right)A(u),
\]
so $A(u(t))$ stays bounded when $M(u_0)<m_d$. At mass $m_d$, the standing wave $e^{it}Q_d$ does not scatter, and its pseudoconformal transform
\[
u(t,x)=|t|^{-d/2}e^{i|x|^2/(4t)-i/t}Q_d(x/t),\qquad t<0,
\]
blows up at time zero \cite{KLR2009,MXZ2010}. Our first application gives a direct $H^1$ classification of finite-time blow-up at this mass.

\begin{theorem}\label{Thm1.2}
Let $d\ge3$ and let $u\in C([0,T);H^1(\R^d;\C))$ solve \eqref{eq1.2}, where $0<T<\infty$ is its maximal forward existence time. If $M(u(0))=m_d$, then there are $x_*,y\in\R^d$, $\theta\in\R$ and $\lambda>0$ such that
\[
u(0,x)=e^{i\theta-i|x-x_*|^2/(4T)}\lambda^{d/2}Q_d\bigl(\lambda(x-y)\bigr).
\]
Consequently, $u$ is obtained from a pseudoconformal standing wave by phase rotations, space and time translations, scaling and Galilean transformations of \eqref{eq1.2}.
\end{theorem}

This is the Hartree counterpart of Merle's theorem \cite{Merle1993}. For $d=4$ it is the classification of Miao, Xu and Zhao \cite[Theorem~1.1]{MXZ2010}, and for radial data with $d\ge5$ it also follows from Theorem~\ref{Thm1.1} and the conditional result of Li and Zhang \cite[Theorem~1.4]{LiZhang2009}. The proof combines Theorem~\ref{Thm1.1} with the quantitative compactness of Corollary~\ref{Cor1.5} below and the localized variance argument of \cite{Banica2004,HmidiKeraani2005,MXZ2010}.

For global solutions, we say that $u$ scatters in $L^2$ as $t\to\pm\infty$ if there is $u_\pm\in L^2(\R^d)$ such that $\|u(t)-e^{it\Delta}u_\pm\|_2\to0$.

\begin{corollary}\label{Cor1.3}
Let $d\ge3$ and let $u\in C(\R;H^1_{\mathrm{rad}}(\R^d;\C))$ be a global solution of \eqref{eq1.2} with $M(u)=m_d$. If $u$ fails to scatter in $L^2$ in at least one time direction, then there are $\theta\in\R$ and $\lambda>0$ such that
\[
u(t,x)=e^{i\theta}e^{i\lambda^2t}\lambda^{d/2}Q_d(\lambda x)
\qquad\text{for every }t\in\R.
\]
\end{corollary}

Corollary~\ref{Cor1.3} follows from the radial classifications of Li and Zhang \cite[Theorem~1.4]{LiZhang2009} for $d\ge5$ and of Gao and Wang \cite[Theorem~1.3]{GaoWang2017} for $d=3,4$, once their hypotheses on the ground state, including the spectral condition are verified. Theorem~\ref{Thm1.1} verifies these hypotheses in every $d\ge3$. It likewise supplies the uniqueness hypothesis in \cite[Appendix~A]{MMRZ2026}, so their Theorem~A.1 applies in every $d\ge3$ to $L^2$ solutions on $[0,\infty)$ that are almost periodic modulo symmetries, have mass $m_d$ and have frequency scale $N(t)\le1$: along a sequence of nonnegative times, the solution converges strongly in $L^2$ to $Q_d$ after phase rotations, translations, scaling and frequency modulation. The proofs of Theorem~\ref{Thm1.2} and Corollary~\ref{Cor1.3} are given at the end of Section~\ref{Sec5}.

Ma, Miao and Zheng \cite[Theorem~1.2]{MMZ2026} recently classified maximal-lifespan strong $L^2$ solutions of mass $m_d$ with infinite $L^{2(d+2)/d}_{t,x}$ norm, assuming uniqueness and nondegeneracy of the ground state. Up to the symmetries of \eqref{eq1.2}, each such solution is a solitary wave or a pseudoconformal transform of one. Theorem~\ref{Thm1.1}, Lemma~\ref{Lem4.4} and the spectral estimates in the proof of Corollary~\ref{Cor1.3} verify their hypotheses, including the four properties in \cite[Theorem~A.1]{MMZ2026}, in every $d\ge3$. Their classification therefore holds unconditionally. It contains Theorem~\ref{Thm1.2}, since the $L^2$ continuation criterion and persistence of regularity imply that an $H^1$ solution that blows up at a finite time $T$ has infinite $L^{2(d+2)/d}_{t,x}$ norm on $[0,T)\times\R^d$. Its restriction to global radial solutions also gives Corollary~\ref{Cor1.3}. Our proof of Theorem~\ref{Thm1.2} is independent of this classification and uses quantitative stability and localized variance.

Miao, Xu and Zhao \cite{MXZ2009} proved scattering below the ground-state mass for radial $L^2$ data, and Ma, Miao, Rosenzweig and Zheng \cite[Theorem~1.6]{MMRZ2026} removed the radial assumption; that proof does not use uniqueness, whereas their Appendix~A does. The four-dimensional multisoliton and blow-up constructions of G\'omez, Schmid and Wu \cite{GSW2026}, and their extension to soliton clusters by Schmid and Wu \cite{SW2026}, use the full linearized kernel and coercivity near the ground state. The uniqueness, kernel and coercivity results are now available in every $d\ge3$; extending those constructions would require further dynamical analysis.

Theorem~\ref{Thm1.1} also yields a quantitative form of \eqref{eq1.5}. By Lemma~\ref{Lem4.1} and Theorem~\ref{Thm1.1}, its nonzero equality cases are exactly
\[
c\,e^{i\theta}\lambda^{d/2}Q_d(\lambda(\cdot-y)),
\qquad c,\lambda>0,\quad \theta\in\R,\quad y\in\R^d.
\]
Amplitude and dilation are not symmetries of \eqref{eq1.1}, whose frequency is normalized to one: the dilation $Q_\lambda(x)=\lambda^{d/2}Q(\lambda x)$ satisfies
\[
-\Delta Q_\lambda+\lambda^2Q_\lambda
=(|x|^{-2}*Q_\lambda^2)Q_\lambda,
\]
and with $\Lambda Q=dQ/2+x\cdot\nabla Q$, differentiation gives
\[
L_+\Lambda Q=-2Q,\qquad \ip{Q}{\Lambda Q}=0.
\]
Thus $\Lambda Q\in Q^\perp$ is a null vector of $q_+$, although $\Lambda Q\notin\ker L_+$. We therefore fix both amplitude and scale, imposing $M(u)=A(u)=m_d$, and measure the distance to
\[
\mathcal O_d=\{e^{i\theta}Q_d(\cdot-y):\theta\in\R,\ y\in\R^d\}.
\]

The estimate below is a Bianchi--Egnell remainder for the sharp Hartree interpolation inequality. For local interpolation inequalities, quantitative results include the Gagliardo--Nirenberg estimates of Carlen and Figalli \cite{CF2013} and Bonforte, Dolbeault, Nazaret and Simonov \cite{BDNS2025}, and the explicit Sobolev stability bounds of Dolbeault, Esteban, Figalli, Frank and Loss \cite{DEFFL2025}. Weinstein's modulation estimates \cite{Weinstein1985} provide a related spectral approach. In the nonlocal setting, d'Avenia and Squassina \cite[Theorem~1.1]{DAS2014} proved a local $H^1$ energy remainder at fixed mass for the three-dimensional Newton interaction $|x|^{-1}$.

\begin{theorem}\label{Thm1.4}
For every $d\ge3$ there is a constant $c_d>0$ such that every $u\in H^1(\R^d;\C)$ with $M(u)=A(u)=m_d$ satisfies
\begin{equation}\label{eq1.6}
2m_d-\Dcal(u)\ge c_d\dist_{H^1}(u,\mathcal O_d)^2.
\end{equation}
\end{theorem}

On the normalized constraint set the deficit in \eqref{eq1.5} equals $4E(u)$, so Theorem~\ref{Thm1.4} bounds the squared $H^1$ distance to $\mathcal O_d$ by a constant multiple of $E(u)$, for complex and nonradial functions. The passage from constrained coercivity to a global remainder follows Bianchi and Egnell \cite{BianchiEgnell1991}; related local modulation estimates for the mass-critical Hartree equation appear in \cite{GaoWang2017}. In dimension four, the elliptic results of \cite{KLR2009,GSW2026} already suffice for this argument. Theorem~\ref{Thm1.1} supplies them in the remaining dimensions. After fixing the scale, the estimate controls concentrating profiles.

\begin{corollary}\label{Cor1.5}
Let $u\in H^1(\R^d;\C)$ satisfy $M(u)=m_d$, and put
\[
\lambda=\left(\frac{m_d}{A(u)}\right)^{1/2},
\qquad v(x)=\lambda^{d/2}u(\lambda x).
\]
Then
\begin{equation}\label{eq1.7}
\dist_{H^1}(v,\mathcal O_d)^2
\le C_d\frac{E(u)}{A(u)}.
\end{equation}
In particular, if $M(u_n)=m_d$, $A(u_n)\to\infty$ and $\sup_n E(u_n)<\infty$, the normalized sequence $\lambda_n^{d/2}u_n(\lambda_n\cdot)$ approaches $\mathcal O_d$ strongly in $H^1$. For a sequence with $M(u_n)\to m_d$, $A(u_n)\to m_d$ and $E(u_n)\to0$, one also has $\dist_{H^1}(u_n,\mathcal O_d)\to0$.
\end{corollary}

The optimal constant in Theorem~\ref{Thm1.4} is
\begin{equation}\label{eq1.8}
\begin{gathered}
\mathcal N=\{u\in H^1(\R^d;\C):M(u)=A(u)=m_d\},\\
\mathcal R(u)=\frac{2m_d-\Dcal(u)}{\dist_{H^1}(u,\mathcal O_d)^2},
\qquad
c_* =\inf_{u\in\mathcal N\setminus\mathcal O_d}\mathcal R(u).
\end{gathered}
\end{equation}
Using the same distance to the full orbit $\mathcal O_d$, define the radial optimal constant by
\[
c_*^{\mathrm{rad}}=
\inf_{u\in(\mathcal N\cap H^1_{\mathrm{rad}}(\R^d;\C))\setminus\mathcal O_d}
\mathcal R(u).
\]
Near the orbit, $\mathcal R$ is governed by two constrained spectral values. Let
\[
\begin{aligned}
&T_+=\{a\in H^1(\R^d;\R):\ \ip{Q_d}{a}=0,
\quad\ip{\nabla Q_d}{\nabla a}=0,
\ip{a}{\partial_{x_j}Q_d}_{H^1}=0\quad(1\le j\le d)\},\\
&T_-=\{b\in H^1(\R^d;\R):\ip{b}{Q_d}_{H^1}=0\},
\end{aligned}
\]
and set
\[
\mu_\pm=\inf_{h\in T_\pm\setminus\{0\}}
\frac{q_\pm[h]}{\|h\|_{H^1}^2}.
\]
Away from the orbit, minimizing sequences may lose compactness through translation, dispersion of mass and concentration of kinetic energy. Vanishing sequences satisfy $\mathcal R(u_n)\to1/2$ by Lemma~\ref{Lem6.6}, so $1/2$ is the dispersive level. Splitting sequences are described by admissible families: possibly empty, at most countable families of nonzero $v^j\in H^1(\R^d;\C)$ such that
\[
\sum_jM(v^j)\le m_d,\qquad
\sum_jA(v^j)\le m_d,\qquad
\sup_j\mathcal C(v^j)<2m_d,
\]
where
\[
\mathcal C(v)=\sup_{y\in\R^d}
\bigl|\ip{v}{Q_d(\cdot-y)}_{H^1}\bigr|.
\]
The relaxed quotient is
\[
\mathcal R_{\mathrm{rel}}((v^j))
=\frac{2m_d-\sum_j\Dcal(v^j)}
{4m_d-2\sup_j\mathcal C(v^j)}.
\]
For the empty family, all sums and the correlation supremum are zero. The interaction is additive over the profiles, while the correlation with $\mathcal O_d$ enters through a supremum. By Lemma~\ref{Lem5.1}, a single profile $v\in\mathcal N\setminus\mathcal O_d$ satisfies $\mathcal R_{\mathrm{rel}}((v))=\mathcal R(v)$.

\begin{theorem}\label{Thm1.6}
Let $d\ge3$, and let $c_*$ be defined by \eqref{eq1.8}. Then:
\begin{enumerate}
\item[(i)] $0<c_*<\tfrac12$.
\item[(ii)] For every nonzero real $b\in H^1(\R^d)$ with $\ip{b}{Q_d}_{H^1}=0$,
\[
c_*<\frac{2q_-[b]}{\|b\|_{H^1}^2}.
\]
In particular, $c_*<2\mu_-$.
\item[(iii)] The local value is
\[
\lim_{\delta\downarrow0}
\inf_{\substack{u\in\mathcal N\\
0<\dist_{H^1}(u,\mathcal O_d)<\delta}}
\mathcal R(u)=2\min\{\mu_+,\mu_-\}.
\]
\item[(iv)] Every minimizing sequence has a subsequence with one of the following properties. Either its distance to $\mathcal O_d$ tends to zero, which is possible only if $c_*=2\mu_+<2\mu_-$; or it admits a translation profile decomposition with a nonempty admissible family minimizing $\mathcal R_{\mathrm{rel}}$, and
\[
\sum_jM(v^j)=m_d\qquad\text{or}\qquad\sum_jA(v^j)=m_d.
\]
In the second case the sequence converges strongly modulo phase and translation if and only if the family consists of one profile $v$ with $M(v)=A(v)=m_d$. Consequently, if $c_*$ is not attained, then either $c_*=2\mu_+<2\mu_-$, or the relaxed problem has a minimizer with a defect: a family with at least two profiles, or a single profile $v$ with $M(v)<m_d$ or $A(v)<m_d$.
\item[(v)] If $v=e^{i\theta}w$ with $\theta\in\R$ and $w\ge0$, and if $M(v)=m_d$ and $A(v)<m_d$, then the one-profile family $(v)$ is not a relaxed minimizer.
\end{enumerate}
\end{theorem}

For the optimal constant, an exact correlation formula for the distance to $\mathcal O_d$ identifies the local value, and a fourth-order expansion in any imaginary tangent direction gives a strict improvement over its Rayleigh quotient. A concentrated ground state together with a diffuse mass correction puts the constant below $1/2$ in every dimension. Translation profile decomposition leads to the relaxed quotient; separate corrections of the missing mass and gradient norm show that every admissible family is realizable, and amplitude and phase variations then constrain relaxed minimizers.

For the Sobolev inequality, K\"onig proved that the optimal stability constant lies strictly below its local spectral value and is attained \cite{Koenig2023,Koenig2025}. Attainment is also known for the Hardy--Sobolev inequality in a range of attractive Hardy parameters \cite{CGK2025} and for the logarithmic Sobolev inequality \cite{DS2026}. For \eqref{eq1.5}, the two constraints can be lost independently, so a minimizing sequence may leave a defect in mass or in kinetic energy. Theorem~\ref{Thm1.6} describes these possible defects and excludes a single profile of constant phase with full mass and deficient kinetic energy. It also leaves open the local alternative $c_*=2\mu_+<2\mu_-$. Whether $c_*$ is attained remains open.

For radial sequences, compactness in $L^{2d/(d-1)}$ \cite{Strauss1977} reduces the relaxed problem to a single weak profile. Theorem~\ref{Thm6.12} gives $c_*^{\mathrm{rad}}<1/2$. If this constant is not attained, then either $c_*^{\mathrm{rad}}=2\mu_+^{\mathrm{rad}}<2\mu_-^{\mathrm{rad}}$, where the spectral infima are restricted to radial functions, or a radial relaxed minimizer has a defect in mass or kinetic energy. In the latter case at least one constraint is saturated, and a constant-phase profile with full mass cannot have deficient kinetic energy.

We now place these results in the context of earlier work. In the nonlocal setting, Lieb \cite{Lieb1977} proved uniqueness of the minimizing solution for the three-dimensional Newton interaction, Lions \cite{Lions1980} developed the variational theory, and Lenzmann \cite{Lenzmann2009} identified the full linearized kernel and applied it to pseudorelativistic Hartree ground states. Arora and Roudenko \cite[Theorem~A.1]{AroraRoudenko2022} extended Newtonian uniqueness to positive radial $H^1$ solutions in dimensions three through five.

For a general Riesz interaction $I_\alpha(x)=c_{d,\alpha}|x|^{\alpha-d}$, $0<\alpha<d$, Moroz and Van Schaftingen \cite{MorozVS2013} established existence, positivity, radial symmetry and decay of ground states in the optimal range $(d+\alpha)/d<p<(d+\alpha)/(d-2)$, and their survey \cite{MorozVS2017} discusses uniqueness and nondegeneracy of positive ground states in this range and their role in semiclassical constructions. If $\Dcal_{\alpha,p}(u)=\int_{\R^d}(I_\alpha*|u|^p)|u|^p$, the scaling $\Dcal_{\alpha,p}(u_\lambda)=\lambda^{dp-d-\alpha}\Dcal_{\alpha,p}(u)$ shows that $p=2$ is mass-critical exactly when $\alpha=d-2$; the Riesz normalization can be absorbed into the amplitude of $Q$. Xiang \cite{Xiang2016} treated the three-dimensional Newton kernel with $p>2$ close to two; Li \cite{Li2022} obtained a priori bounds and both properties near $(\alpha,p)=(2,2)$ for $d\in\{3,4,5\}$; Luo, Zhang and Xu \cite{LZX2024} treated $p=2$ in dimension three for $\alpha$ close to two. Seok \cite[Theorem~1.5 and Remark~1.6]{Seok2019} proved uniqueness and nondegeneracy near $\alpha=0$ when $2\le p<d/(d-2)$, and near $\alpha=d$ when $2<p<2d/(d-2)$; for $d\ge3$ the first range is nonempty only in dimension three, and the second excludes $p=2$.

For the boson star equation $(\sqrt{-\Delta}+1)Q=(|x|^{-1}*Q^2)Q$ in $\R^3$, Chen, Guo, Wei and Yu \cite{CGWY2026} proved uniqueness and nondegeneracy of ground states; that equation has fractional dispersion and a Newton interaction, whereas \eqref{eq1.1} has the ordinary Laplacian and a fixed inverse-square interaction.

Section~\ref{Sec2} develops the radial family and its a priori bounds. Sections~\ref{Sec3} and \ref{Sec4} prove nondegeneracy, uniqueness and the full kernel characterization. Section~\ref{Sec5} proves the stability estimate and the results on threshold dynamics. Section~\ref{Sec6} treats the optimal stability constant.

\section{The radial family and its compactness}\label{Sec2}

We embed the radial equation in a family with real dimension parameter and prove compactness of its positive decreasing solutions. Fix $\nu_*>4$ and let $3\le\nu\le \nu_*$. Constants denoted by $C_{\nu_*}$ depend only on $\nu_*$. Set
\[
\omega_\nu=\frac{2\pi^{\nu/2}}{\Gamma(\nu/2)},\qquad
\kappa_\nu=\frac{2\pi^{(\nu-1)/2}}{\Gamma((\nu-1)/2)},\qquad
\gamma_\nu=\frac{(\nu-1)(\nu-3)}4.
\]
For $r,s>0$ define
\begin{equation}\label{eq2.1}
\begin{aligned}
a_\nu(r,s)&=\kappa_\nu\int_{-1}^1
\frac{(1-t^2)^{(\nu-3)/2}}{r^2+s^2-2rst}\dd t,\\
a_{\nu,1}(r,s)&=\kappa_\nu\int_{-1}^1
\frac{t(1-t^2)^{(\nu-3)/2}}{r^2+s^2-2rst}\dd t,\\
k_\nu(r,s)&=\kappa_\nu\int_{-1}^1
\frac{t(1-t^2)^{(\nu-3)/2}}{(r^2+s^2-2rst)^2}\dd t.
\end{aligned}
\end{equation}
The last integral is defined away from $r=s$. Diagonal values of locally integrable kernels do not affect the associated integral operators. Write
\[
\mathcal V_\nu[z](r)=\int_0^\infty a_\nu(r,s)z(s)^2\dd s,
\qquad \Ncal_\nu(z)=\mathcal V_\nu[z]z.
\]
For the bilinear expression obtained by replacing $z^2$ with $z_1z_2$, we write $\mathcal V_\nu[z_1,z_2]$.

On the fixed real Hilbert space $\Hh=H_0^1(0,\infty)$, with norm $\|z\|_\Hh^2=\int(z'^2+z^2)$, consider
\begin{equation}\label{eq2.2}
\Fcal(\nu,z)=A_\nu z-\Ncal_\nu(z)=0,
\qquad
A_\nu=-\frac{\mathrm d^2}{\mathrm dr^2}+1+\frac{\gamma_\nu}{r^2}.
\end{equation}
The operator $A_\nu$ is interpreted through its form from $\Hh$ to $\Hh^*$. The half-line Hardy inequality
\begin{equation}\label{eq2.3}
\int_0^\infty\frac{z^2}{r^2}\dd r\le4\int_0^\infty z'^2\dd r
\end{equation}
shows that $A_\nu$ is uniformly bounded and coercive on the stated parameter interval.

If $Q(r)=r^{-(\nu-1)/2}z(r)$, direct conjugation gives
\begin{equation}\label{eq2.4}
-Q''-\frac{\nu-1}{r}Q'+Q=V_QQ,
\qquad
V_Q(r)=\int_0^\infty a_\nu(r,s)Q(s)^2s^{\nu-1}\dd s.
\end{equation}
Moreover,
\begin{equation*}
\int_0^\infty(Q'^2+Q^2)r^{\nu-1}\dd r
=\int_0^\infty\left(z'^2+z^2+\frac{\gamma_\nu}{r^2}z^2\right)\dd r.
\end{equation*}
The norm identity holds first for functions supported in compact annuli and extends by completion. At integer $\nu=d$, the Euclidean quantities are
\[
M(Q)=\omega_d\|z\|_2^2,\qquad
A(Q)=\omega_d\int_0^\infty
\left(z'^2+\frac{\gamma_d}{r^2}z^2\right)\dd r.
\]
Equation \eqref{eq2.4} is then the radial form of \eqref{eq1.1}. For noninteger $\nu$, \eqref{eq2.2} defines the half-line equation used below.

\begin{lemma}\label{Lem2.1}
The kernels in \eqref{eq2.1} have the following properties.

\smallskip
\noindent\textup{(i)} The kernel $a_\nu$ is symmetric, positive and homogeneous of degree $-2$. For $r\ne s$ and $j=0,1$,
\begin{equation}\label{eq2.5}
|\partial_\nu^j a_\nu(r,s)|
\le\frac{C_{\nu_*,j}}{(r+s)^2}
\left(1+\left|\log\frac{|r-s|}{r+s}\right|^{j+1}\right).
\end{equation}
On $s<r/2$ or $s>2r$, the logarithmic factor may be omitted.

\smallskip
\noindent\textup{(ii)} For $r,s>0$ with $r\ne s$, $a_{\nu,1}(r,s)>0$. If $0<\theta_1<\theta_2<\infty$, then
\begin{equation}\label{eq2.6}
a_{\nu,1}(r,s)\ge c_{\nu_*,\theta_1,\theta_2}r^{-2}
\qquad\text{when }\theta_1\le s/r\le\theta_2.
\end{equation}
In the distributional sense on $(0,\infty)^2$,
\begin{equation}\label{eq2.7}
\partial_r a_\nu(r,s)
=-\left(\partial_s+\frac{\nu-1}{s}\right)a_{\nu,1}(r,s).
\end{equation}

\smallskip
\noindent\textup{(iii)} Let $B_\nu=\partial_r\partial_s a_\nu$. Away from $r=s$,
\begin{equation}\label{eq2.8}
B_\nu(r,s)=2(\nu-4)k_\nu(r,s),\qquad k_\nu(r,s)>0.
\end{equation}
For $\nu>4$, this is an ordinary locally integrable kernel. For $\nu=4$,
\[
a_4(r,s)=\omega_4\max\{r,s\}^{-2},
\qquad B_4(r,s)=2\omega_4r^{-3}\delta(r-s).
\]
For $3\le\nu<4$, \eqref{eq2.8} describes only the off-diagonal part of the distribution.
\end{lemma}

\begin{proof}
Symmetry and homogeneity are immediate. Put
\[
b=r^2+s^2-2rst=(r-s)^2+2rs(1-t).
\]
For $s<r/2$ or $s>2r$, the denominator is comparable to $\max(r,s)^2$. Differentiation in $\nu$ introduces derivatives of $\kappa_\nu$ and powers of $\log(1-t^2)$, whose weighted angular integrals are uniformly bounded. For $r/2\le s\le2r$, divide by $(r+s)^2$. Near $t=1$, the integral is bounded by a constant times
\[
\int_0^1\frac{1+|\log u|^j}{\bigl(|r-s|/(r+s)\bigr)^2+u}\dd u
\le C_j\left(1+\left|\log\frac{|r-s|}{r+s}\right|^{j+1}\right).
\]
The other endpoint is integrable. This proves \eqref{eq2.5}.

Pairing the values at $t$ and $-t$ proves positivity of the two angular integrals with numerator $t$. Restriction to $1/4<t<1/2$ gives \eqref{eq2.6} uniformly for $s/r$ in a compact subset of $(0,\infty)$.

The identity underlying \eqref{eq2.7} is
\[
\partial_r(b^{-1})+\partial_s(tb^{-1})
+\frac{\nu-1}{s}tb^{-1}
=-\frac{1}{s(1-t^2)^{(\nu-3)/2}}
\partial_t\left(\frac{(1-t^2)^{(\nu-1)/2}}{b}\right).
\]
Replace $b$ by $b+\varepsilon^2$ in the identity and integrate in $t$. The angular boundary terms vanish. Local integrability of $a_\nu$ and $a_{\nu,1}$ then allows $\varepsilon\downarrow0$ in distributions, proving \eqref{eq2.7}.

For the mixed derivative the corresponding identity is
\[
\partial_r\partial_s(b^{-1})
-2(\nu-4)\frac{t}{b^2}
=\frac{2}{(1-t^2)^{(\nu-3)/2}}
\partial_t\left(\frac{(1-t^2)^{(\nu-1)/2}}{b^2}\right).
\]
For $r\ne s$ the angular boundary terms vanish, giving \eqref{eq2.8}. On a compact annulus the singular behavior of $k_\nu$ is bounded by
\[
C_\nu\begin{cases}
1+|r-s|^{\nu-5},&3\le\nu<5,\\
1+|\log|r-s||,&\nu=5,\\
1,&\nu>5.
\end{cases}
\]
The bound is locally integrable for $\nu>4$. If $4<\nu<5$, the nonsmooth part of $a_\nu$ has order $|r-s|^{\nu-3}$, and its first derivative tends to zero at the diagonal. At $\nu=5$, that derivative has order $|r-s|\,|\log|r-s||$. For $\nu>5$, differentiation under the angular integral gives continuity. Thus the first derivative has no jump, and \eqref{eq2.8} represents the full mixed derivative for $\nu>4$. At $\nu=4$, direct integration gives the shell formula and the stated diagonal measure. For $3\le\nu<4$, the mixed derivative is retained in the distributional sense. In particular,
\[
a_3(r,s)=\frac{2\pi}{rs}\log\frac{r+s}{|r-s|}.
\]
\end{proof}

\begin{lemma}\label{Lem2.2}
The map $\Ncal_\nu:\Hh\to L^2(0,\infty)$ is continuously differentiable, with
\begin{equation}\label{eq2.9}
D\Ncal_\nu(z)y=\mathcal V_\nu[z]y
+2\mathcal V_\nu[z,y]z.
\end{equation}
Both the map and its derivative depend continuously on $\nu$, with the derivative continuous in operator norm. In addition,
\begin{equation}\label{eq2.10}
\|\mathcal V_\nu[z]\|_\infty\le C_{\nu_*}\|z'\|_2^2,
\qquad
|\mathcal V_\nu[z](r)|
\le C_{\nu_*}\frac{1+\log r}{r^2}\|z\|_\Hh^2\quad(r\ge2).
\end{equation}
If $\nu_n\to\nu$ and $z_n\rightharpoonup z$ in $\Hh$, then
\begin{equation}\label{eq2.11}
\Ncal_{\nu_n}(z_n)\longrightarrow\Ncal_\nu(z)
\quad\text{in }L^2(0,\infty).
\end{equation}
The derivative $D_z\Fcal(\nu,z):\Hh\to\Hh^*$ is Fredholm of index zero.
\end{lemma}

\begin{proof}
For $z\in\Hh$,
\begin{equation}\label{eq2.12}
|z(s)|^2\le s\|z'\|_2^2,
\qquad \int_0^\infty z(s)^2s^{-2}\dd s\le4\|z'\|_2^2.
\end{equation}
Split the potential integral into $s<r/2$, $r/2<s<2r$ and $s>2r$. On the first region, combine $a_\nu\le C_{\nu_*}r^{-2}$ with the pointwise bound in \eqref{eq2.12}. On the third, use $a_\nu\le C_{\nu_*}s^{-2}$ and the Hardy bound. On the middle region, apply \eqref{eq2.5}, $|z(s)|^2\le2r\|z'\|_2^2$ and integrability of $1+|\log|1-s/r||$. This proves the first estimate in \eqref{eq2.10}. The same argument applies to the first parameter derivative. Cauchy--Schwarz with the absolute kernel bound gives the bilinear estimates for $\mathcal V_\nu[z_1,z_2]$.

For $r\ge2$, the regions $s<r/2$ and $s>2r$ contribute at most $C_{\nu_*}r^{-2}\|z\|_2^2$. In the middle region with $|s-r|\ge1$, the kernel is bounded by $C_{\nu_*}(1+\log r)/r^2$. For $|s-r|<1$, use $\|z\|_\infty\le C\|z\|_\Hh$ and integrate the logarithmic singularity. This gives the second estimate in \eqref{eq2.10}. Expansion of the cubic map $\Ncal_\nu$ and the bilinear estimates give \eqref{eq2.9}. Differentiation of the angular kernel in $\nu$ also gives
\[
\|\mathcal V_\nu[z_1,z_2]-\mathcal V_\mu[z_1,z_2]\|_\infty
\le C_{\nu_*}|\nu-\mu|\|z_1'\|_2\|z_2'\|_2.
\]
Together with the bilinear estimates, this proves continuity in the profile and parameter, including operator-norm continuity of the derivative.

Under the hypotheses of \eqref{eq2.11}, local compactness gives a subsequence converging uniformly to $z$ on compact subintervals of $(0,\infty)$. The potentials converge uniformly on those intervals as well. To see this, use \eqref{eq2.12} for small $s$, the separated-radius estimate for large $s$, and the integrable logarithmic majorant in \eqref{eq2.5} for the remaining region. The same estimate controls the dependence on $\nu$. Hence the nonlinearities converge locally in $L^2$. Near zero,
\[
\|\Ncal_{\nu_n}(z_n)\|_{L^2(0,\varepsilon)}
\le C_{\nu_*}\|z_n\|_{L^2(0,\varepsilon)}
\le C_{\nu_*}\varepsilon.
\]
At infinity, the second estimate in \eqref{eq2.10} gives uniformly small $L^2$ tails. Thus \eqref{eq2.11} holds along the extracted subsequence. Applying the same argument to every subsequence proves convergence of the full sequence.

The compactness of the derivative follows directly from the cubic identity
\[
D\Ncal_\nu(z)y
=\frac12\bigl(\Ncal_\nu(z+y)-\Ncal_\nu(z-y)\bigr)-\Ncal_\nu(y).
\]
Each term on the right has relatively compact image on bounded subsets of $\Hh$. Hence the linear map $D\Ncal_\nu(z):\Hh\to L^2$ is compact.

Finally, the form of $A_\nu$ is bounded and satisfies
\[
\ip{A_\nu y}{y}\ge\|y\|_\Hh^2,\qquad
|\ip{(A_\nu-A_\mu)y}{h}|
\le4|\gamma_\nu-\gamma_\mu|\|y'\|_2\|h'\|_2.
\]
The Lax--Milgram theorem gives an isomorphism $A_\nu:\Hh\to\Hh^*$ with uniformly bounded inverse. The resolvent identity gives continuity of $A_\nu^{-1}$ in $\nu$. Since $D_z\Fcal=A_\nu-D\Ncal_\nu(z)$ is a compact perturbation of $A_\nu$, it is Fredholm of index zero.
\end{proof}

Let $\mathcal P_\nu$ be the set of nonzero solutions of \eqref{eq2.2} for which
\[
Q(r)=r^{-(\nu-1)/2}z(r)
\]
has a nonnegative, nonincreasing representative. Lemma~\ref{Lem2.3} proves that this representative is positive and strictly decreasing. For such a solution, put $V=\mathcal V_\nu[z]$ and define
\begin{equation}\label{eq2.13}
\Lcal_\nu y=A_\nu y-Vy-2z\mathcal V_\nu[z,y],
\qquad q_\nu[y_1,y_2]=\ip{\Lcal_\nu y_1}{y_2}.
\end{equation}
Unweighted pairings on the half-line use $\mathrm dr$, and $q_\nu[y]=q_\nu[y,y]$. The Morse index is the maximal dimension of a subspace on which $q_\nu$ is negative definite. Lemma~\ref{Lem2.2} implies that it is finite.

\begin{lemma}\label{Lem2.3}
Let $z\in\mathcal P_\nu$. Then $Q(0)>0$, $Q'(0)=0$, $Q'(r)<0$ and $V'(r)<0$ for $r>0$. At the origin,
\[
Q(r)=Q(0)+O(r^2),\qquad Q'(r)=O(r).
\]
At infinity, for some $a>0$,
\begin{equation}\label{eq2.14}
\begin{gathered}
z(r)=ae^{-r}\left(1+O\left(\frac{1+\log r}{r}\right)\right),\qquad
\frac{z'(r)}{z(r)}=-1+O\left(\frac{1+\log r}{r^2}\right).
\end{gathered}
\end{equation}
If
\[
\rho=z^2,\qquad Z=\frac12z+rz',\qquad
\eta=\frac Zz=\frac12+\frac{r\rho'}{2\rho},
\]
then $Z\in\Hh$ and, at infinity,
\[
\eta(r)=\frac12-r+O\left(\frac{1+\log r}{r}\right),\qquad
\eta'(r)=-1+O\left(\frac{1+\log r}{r}\right),\qquad
\int_r^\infty\rho(s)\dd s\le C\rho(r).
\]
Finally,
\begin{equation}\label{eq2.15}
\int_0^\infty\left(z'^2+\frac{\gamma_\nu}{r^2}z^2\right)\dd r
=\|z\|_2^2,\qquad
\int_0^\infty Vz^2\dd r=2\|z\|_2^2,
\end{equation}
and
\begin{equation}\label{eq2.16}
\Lcal_\nu Z=-2z,\qquad
\Lcal_\nu z=-2Vz,\qquad \ip{z}{Z}=0.
\end{equation}
\end{lemma}

\begin{proof}
By Lemma~\ref{Lem2.2}, $V$ is bounded and continuous at positive radii. To prove continuity at zero, split its integral into $s<\varepsilon$ and $s\ge\varepsilon$. The local version of \eqref{eq2.12} makes the first part uniformly small for small $r$, while dominated convergence applies to the second.

On a sufficiently small interval the equation for $Q$ has a regular solution $\phi$ with $\phi(0)=1$, given by
\begin{equation}\label{eq2.17}
\phi(r)=1+\int_0^r t^{1-\nu}
\int_0^t (1-V(s))\phi(s)s^{\nu-1}\dd s\dd t.
\end{equation}
The integral operator is a contraction on a small interval, and $\phi=1+O(r^2)$, $\phi'=O(r)$. A second independent solution is
\[
\phi(r)\int_{r_0}^r\frac{t^{1-\nu}}{\phi(t)^2}\dd t,
\]
whose leading term is a nonzero multiple of $r^{2-\nu}$. After multiplication by $r^{(\nu-1)/2}$, this solution has nonzero trace at $\nu=3$ and a derivative outside $L^2$ for $\nu>3$, so it is excluded by $z\in\Hh$. Thus $Q$ extends continuously to zero with the asserted expansion. Since $Q$ is nonnegative, nonincreasing and nonzero, $Q(0)=\sup Q>0$. The strong maximum principle on compact annuli implies $Q>0$ on $(0,\infty)$.

By \eqref{eq2.10}, $V(r)\to0$. Choose $R$ so that $V\le1/4$ on $[R,\infty)$. The function $Q(R)e^{-(r-R)/2}$ is a supersolution of \eqref{eq2.4}, and the exterior maximum principle gives
\begin{equation*}
Q(r)\le Q(R)e^{-(r-R)/2}\qquad(r\ge R).
\end{equation*}
Testing the positive part of the difference after cutting off at large radii proves the comparison. Estimates for the equation on unit intervals give the same exponential decay, up to polynomial factors, for the first derivatives.

Apply \eqref{eq2.7} to $Q^2$. Integration by parts, first with compact cutoffs, yields
\begin{equation}\label{eq2.18}
-V'(r)=\int_0^\infty
a_{\nu,1}(r,s)(-(Q^2)'(s))s^{\nu-1}\dd s.
\end{equation}
The estimate $(Q^2)'=O(s)$ at zero and exponential decay at infinity allow the cutoffs to be removed. Local integrability of the angular kernel, or its positive regularization from Lemma~\ref{Lem2.1}, justifies integration across $s=r$. Since $Q^2$ is nonincreasing and nonconstant, the positive kernel in \eqref{eq2.18} gives $V'<0$.

The radial equation implies $V(0)>1$. Otherwise $(r^{\nu-1}Q')'=(1-V)Qr^{\nu-1}\ge0$, contradicting positivity and decay of $Q$. To obtain strict decrease at every positive radius, put $y=-Q'\ge0$. Differentiation on an annulus gives
\[
-y''-\frac{\nu-1}{r}y'
+\left(1-V+\frac{\nu-1}{r^2}\right)y=-V'Q>0.
\]
At an interior zero of $y$, the left side is nonpositive. Hence $Q'<0$ on $(0,\infty)$.

For large $r$, write the equation for $z$ as $z''=(1+p)z$, where
\[
p(r)=\gamma_\nu r^{-2}-V(r).
\]
Since $V'<0$, \eqref{eq2.10} implies
\[
|p(r)|+\int_r^\infty|p'(s)|\dd s
\le 2\bigl(\gamma_\nu r^{-2}+V(r)\bigr)
\le C\frac{1+\log r}{r^2}.
\]
We have $z''>0$ for large $r$ and $z,z'\to0$ by the preceding decay estimates, so $z'<0$ there. Integrating
\[
\bigl(z'^2-(1+p)z^2\bigr)'=-p'z^2
\]
from $r$ to infinity gives
\[
\left(\frac{z'}z\right)^2
=1+p(r)+\frac1{z(r)^2}\int_r^\infty p'(s)z(s)^2\dd s
=1+O\left(\frac{1+\log r}{r^2}\right),
\]
where we used $z(s)\le z(r)$ for $s\ge r$. Thus, with $u=z'/z$,
\[
u=-1+O\left(\frac{1+\log r}{r^2}\right),\qquad
u'=1+p-u^2=O\left(\frac{1+\log r}{r^2}\right).
\]
The function $u+1$ is integrable at infinity. Integrating $(\log(e^rz))'=u+1$ shows that $e^rz$ has a finite positive limit $a$ and proves \eqref{eq2.14}. The asserted bounds for $\eta$ and $\eta'$ follow from $\eta=1/2+ru$. Since $u\le-1/2$ eventually, $\rho(s)\le\rho(r)e^{-(s-r)}$ for $s\ge r$ sufficiently large, which proves the bound for $\int_r^\infty\rho$. The origin estimates and the equation for $z$ also imply $Z\in\Hh$.

Testing \eqref{eq2.2} against $z$ equates the sum of the kinetic and mass terms with the interaction. Under $z_\lambda(r)=\lambda^{1/2}z(\lambda r)$, the kinetic term and the interaction both acquire the factor $\lambda^2$, while the mass is unchanged. Differentiating the action along this curve shows that the interaction equals twice the kinetic term. These identities give \eqref{eq2.15}, and the endpoint estimates justify the differentiation. The scaled equation is
\[
\left(-\frac{\mathrm d^2}{\mathrm dr^2}
+\frac{\gamma_\nu}{r^2}+\lambda^2\right)z_\lambda
=\Ncal_\nu(z_\lambda).
\]
Differentiation at $\lambda=1$ gives the first identity in \eqref{eq2.16}. The second follows by substituting $z$ in \eqref{eq2.13}, and integration of $(rz^2)'$ gives the last.
\end{proof}

\begin{proposition}\label{Prop2.4}
There are constants $0<c_{\nu_*}<C_{\nu_*}<\infty$ such that every $z\in\mathcal P_\nu$, with $3\le\nu\le \nu_*$, satisfies
\begin{equation}\label{eq2.19}
c_{\nu_*}\le\|z\|_\Hh\le C_{\nu_*}.
\end{equation}
\end{proposition}

\begin{proof}
Constants in this proof depend only on $\nu_*$ and the fixed auxiliary parameters. Put
\[
m(r)=\int_0^r Q(s)^2s^{\nu-1}\dd s.
\]
The local ground-state transform associated with the positive solution $Q$ gives
\[
\int_0^\infty V\phi^2r^{\nu-1}\dd r
\le\int_0^\infty(\phi'^2+\phi^2)r^{\nu-1}\dd r
\quad\text{for }\phi\in C_c^\infty(0,\infty).
\]
Choose a fixed nonzero $\phi$ supported in $(1/4,3/4)$. Since $V$ is decreasing, this yields $V(1)\le C_{\nu_*}$.

The lower bound $a_\nu(r,s)\ge\omega_\nu(r+s)^{-2}$ and the three-region upper bound, using monotonicity of $Q^2$ in the comparable-radius region, imply
\begin{equation}\label{eq2.20}
V(r)\asymp_{\nu_*} r^{-2}m(r)
+\int_r^\infty Q(s)^2s^{\nu-3}\dd s.
\end{equation}
In the region $r/2<s<2r$,
\[
\int_{r/2}^{2r}a_\nu(r,s)Q(s)^2s^{\nu-1}\dd s
\le C_{\nu_*} r^{\nu-2}Q(r/2)^2
\le C_{\nu_*}r^{-2}m(r).
\]
The last inequality follows by integrating $Q(s)^2\ge Q(r/2)^2$ on $(0,r/2)$. On the remaining regions the kernel is comparable to the inverse square of the larger radius, giving \eqref{eq2.20}.

The bound for $V(1)$ also implies
\begin{equation}\label{eq2.21}
m(r)\le C_{\nu_*}r^2,\qquad
Q(r)^2\le C_{\nu_*}r^{2-\nu}\qquad(r\ge1).
\end{equation}
Indeed, $V(1)\ge c_{\nu_*}(1+r)^{-2}m(r)$, and monotonicity gives $m(r)\ge r^\nu Q(r)^2/\nu$.

Let $R>0$ be the unique radius for which $V(R)=1/4$. It exists because $V(0)>1$, $V'<0$ and $V(\infty)=0$. By exterior comparison,
\begin{equation}\label{eq2.22}
Q(r)\le Q(R)e^{-(r-R)/2}\qquad(r\ge R).
\end{equation}
If $R\ge1$, the tail in \eqref{eq2.20}, evaluated at $R$, satisfies
\[
\int_R^\infty Q(s)^2s^{\nu-3}\dd s
\le C_{\nu_*}Q(R)^2R^{\nu-3}
\le C_{\nu_*}m(R)R^{-3}.
\]
Thus
\begin{equation}\label{eq2.23}
c_{\nu_*}R^2\le m(R)\le C_{\nu_*}R^2.
\end{equation}

Suppose, for a contradiction, that $R$ is unbounded as $3\le\nu\le \nu_*$ and $z\in\mathcal P_\nu$. Choose $L>2$ sufficiently large that \eqref{eq2.21} and \eqref{eq2.23} give
\[
m(R)-m(R/L)\ge c_{\nu_*}R^2
\]
whenever $R\ge2L$. Since $Q(s)\le Q(R/L)$ on this annulus, its mass gives a lower bound for $Q(R/L)$. Combining this with monotonicity and \eqref{eq2.21}, we obtain
\begin{equation}\label{eq2.24}
c_{\nu_*} R^{(2-\nu)/2}\le Q(r)\le C_{\nu_*} R^{(2-\nu)/2}
\qquad\left(\frac{R}{2L}\le r\le\frac RL\right).
\end{equation}
The constants depend only on $\nu_*$ and the chosen $L$.

For fixed $0<\theta_1<1<\theta_2$, restrict \eqref{eq2.18} to $\theta_1r\le s\le\theta_2r$. There, \eqref{eq2.6} gives $a_{\nu,1}\ge c_{\nu_*,\theta_1,\theta_2}r^{-2}$, and $s^{\nu-1}\ge c_{\nu_*,\theta_1,\theta_2}r^{\nu-1}$. Integrating $-(Q^2)'\ge0$ yields
\begin{equation*}
-V'(r)\ge c_{\nu_*,\theta_1,\theta_2}r^{\nu-3}
\bigl(Q(\theta_1r)^2-Q(\theta_2r)^2\bigr).
\end{equation*}
Take $\theta_1=1/2$ and $\theta_2=4L$. For $r\in I_R=[R/(2L),R/L]$ we have $\theta_2r\ge2R$. Estimates \eqref{eq2.21}, \eqref{eq2.22} and \eqref{eq2.24} show, for large $R$, that
\begin{equation}\label{eq2.25}
-V'(r)\ge c_{\nu_*}/R\qquad(r\in I_R).
\end{equation}
Hence $V$ decreases by a fixed positive amount across any fixed proportion of $I_R$.

If $V(3R/(4L))\ge1$, then \eqref{eq2.25} gives $V\ge1+\varepsilon_{\nu_*}$ on the leftmost quarter of $I_R$. If $V(3R/(4L))\le1$, it gives $V\le1-\varepsilon_{\nu_*}$ on the rightmost quarter. Put $y=-Q'\ge0$. The radial equation reads
\begin{equation}\label{eq2.26}
(r^{\nu-1}y)'=(V-1)Qr^{\nu-1}.
\end{equation}
Denote the selected quarter by $[a_R,b_R]$. In the first case, \eqref{eq2.26} gives
\[
r^{\nu-1}y(r)\ge c_{\nu_*} R^{(2-\nu)/2}\int_{a_R}^r s^{\nu-1}\dd s
\qquad(a_R\le r\le b_R).
\]
In the second case, integration backwards gives
\[
r^{\nu-1}y(r)\ge c_{\nu_*} R^{(2-\nu)/2}\int_r^{b_R}s^{\nu-1}\dd s.
\]
Since $y\ge0$ at the endpoint of integration and all radii are comparable to $R$, the first estimate on the right half of the selected quarter, or the second on its left half, gives
\[
y(r)\ge c_{\nu_*} R\,R^{(2-\nu)/2}
\]
on a subinterval of length at least $c_{\nu_*}R$. A second integration gives a decrease of $Q$ of at least $c_{\nu_*}R^2R^{(2-\nu)/2}$, contradicting \eqref{eq2.24} for large $R$. Hence
\begin{equation}\label{eq2.27}
R\le C_{\nu_*}.
\end{equation}

The case $R<1$ already satisfies \eqref{eq2.27} after increasing $C_{\nu_*}$ to at least one. Choose $R_{\nu_*}>1$ larger than this uniform bound. Then $V\le1/4$ on $[R_{\nu_*},\infty)$ in both cases. Estimate \eqref{eq2.21} bounds $m(R_{\nu_*})$ and $Q(R_{\nu_*})$, while exterior comparison from $R_{\nu_*}$ bounds the remaining mass. Thus $\|z\|_2^2\le C_{\nu_*}$. The identities \eqref{eq2.15} and $\gamma_\nu\ge0$ give the upper bound in \eqref{eq2.19}. For the lower bound, test \eqref{eq2.2} against $z$ and use \eqref{eq2.10} to obtain
\[
\|z\|_\Hh^2\le\ip{A_\nu z}{z}
=\int_0^\infty Vz^2\dd r\le C_{\nu_*}\|z\|_\Hh^4.
\]
Since $z\ne0$, the lower bound follows.
\end{proof}

\begin{lemma}\label{Lem2.5}
Suppose $\nu_n\to\nu\in[3,\nu_*]$ and $z_n\in\mathcal P_{\nu_n}$. Then a subsequence converges strongly in $\Hh$ to some $z\in\mathcal P_\nu$.
\end{lemma}

\begin{proof}
By Proposition~\ref{Prop2.4}, a subsequence converges weakly to $z$ in $\Hh$. Lemma~\ref{Lem2.2} gives strong $L^2$ convergence of the nonlinearities. The Hardy estimate \eqref{eq2.3} gives operator-norm convergence $A_{\nu_n}\to A_\nu$ from $\Hh$ to $\Hh^*$. Their inverses are uniformly bounded, so
\[
z_n=A_{\nu_n}^{-1}\Ncal_{\nu_n}(z_n)
\longrightarrow A_\nu^{-1}\Ncal_\nu(z)=z
\quad\text{in }\Hh.
\]
The lower bound in \eqref{eq2.19} excludes $z=0$. On compact annuli, the corresponding $Q_n$ converge uniformly to $Q=r^{-(\nu-1)/2}z$, so nonnegativity and monotonicity pass to the limit. Lemma~\ref{Lem2.3} gives strict positivity and strict decrease.
\end{proof}

\Needspace{7\baselineskip}
\section{Radial nondegeneracy}\label{Sec3}

Let $z\in\mathcal P_\nu$. We use the linearization $\Lcal_\nu$ and its form $q_\nu$ from \eqref{eq2.13}. We first deduce nonnegativity on the mass tangent space from the Morse index.

\begin{lemma}\label{Lem3.1}
If $\ind(\Lcal_\nu)=1$, then
\begin{equation*}
q_\nu[y]\ge0\qquad\text{whenever }y\in\Hh,\quad\ip{z}{y}=0.
\end{equation*}
\end{lemma}

\begin{proof}
Suppose $y\perp z$ and $q_\nu[y]<0$. Set
\[
h=z-\frac{q_\nu[z,y]}{q_\nu[y]}y.
\]
Then $q_\nu[h,y]=q_\nu[Z,y]=0$. By \eqref{eq2.16},
\[
q_\nu[Z]=0,\qquad q_\nu[Z,h]=-2\|z\|_2^2\ne0.
\]
The Gram determinant of the form on $\spanop\{Z,h\}$ is negative. This subspace contains a negative direction that is $q_\nu$-orthogonal to $y$. Together with $y$, it spans a two-dimensional negative subspace, contradicting the Morse index.
\end{proof}

For $y\perp z$ put
\begin{equation*}
\rho=z^2,\qquad F_y(r)=\int_0^r z(s)y(s)\dd s,
\qquad w(r)=\frac12r\rho(r),\qquad
\mathcal A=-\frac{\mathrm d}{\mathrm dr}\left(\rho^{-1}\frac{\mathrm d}{\mathrm dr}\right).
\end{equation*}
Then $F_y(0)=F_y(\infty)=0$ and $y/z=F_y'/\rho$. The cumulative function associated with the dilation derivative satisfies
\begin{equation*}
w=F_Z>0,\qquad \frac{w'}{\rho}=\eta,
\qquad \mathcal A w=-\eta'.
\end{equation*}
For $F\in C_c^\infty(0,\infty)$, set $y_F=F'/z$. Then $y_F\in\Hh$ and $y_F\perp z$.

\begin{lemma}\label{Lem3.2}
For $F\in C_c^\infty(0,\infty)$,
\begin{equation}\label{eq3.1}
q_\nu[y_F]=\int_0^\infty\rho(\mathcal A F)^2\dd r
-2\ip{B_\nu}{F\otimes F},
\end{equation}
where the distributional pairing is
\[
\ip{B_\nu}{F\otimes G}
=\iint_{(0,\infty)^2}a_\nu(r,s)F'(r)G'(s)\dd r\dd s.
\]
Furthermore, every $y\in\ker\Lcal_\nu$ satisfies
\begin{equation*}
\ip{z}{y}=0,\qquad \ip{Vz}{y}=0.
\end{equation*}
\end{lemma}

\begin{proof}
The local equation $(A_\nu-V)z=0$ gives the ground-state identity
\[
\ip{(A_\nu-V)y_F}{y_F}
=\int_0^\infty\rho\left|\left(\frac{y_F}{z}\right)'\right|^2\dd r.
\]
Substituting $zy_F=F'$ into the nonlocal term gives \eqref{eq3.1}. Local integrability of $a_\nu$ justifies the pairing on compact annuli. If $\Lcal_\nu y=0$, test against $Z$ and $z$, which belong to $\Hh$ by Lemma~\ref{Lem2.3}. Symmetry of the form and \eqref{eq2.16} give the two orthogonality relations.
\end{proof}

\Needspace{9\baselineskip}
\begin{proposition}\label{Prop3.3}
For $\nu>4$, every $z\in\mathcal P_\nu$ has
\[
\ker\Lcal_\nu=\{0\}.
\]
Its radial Morse index is one. More precisely, for $y\in\Hh$ with $y\perp z$,
\begin{equation}\label{eq3.2}
q_\nu[y]\ge
\iint_{(0,\infty)^2}B_\nu(r,s)w(r)w(s)
\left|\frac{F_y(r)}{w(r)}-\frac{F_y(s)}{w(s)}\right|^2\dd r\dd s.
\end{equation}
\end{proposition}

\begin{proof}
In this parameter range $B_\nu$ is an ordinary positive kernel. Set
\[
U=\rho \mathcal A w=-\rho\eta',\qquad
b(r)=2\int_0^\infty B_\nu(r,s)w(s)\dd s.
\]
For $s>2r$, $B_\nu(r,s)=O_\nu(rs^{-5})$, and symmetry gives the estimate for $s<r/2$. On comparable radii, use homogeneity and local integrability from Lemma~\ref{Lem2.1}. Since $w(s)=O(s^\nu)$ at zero and decays exponentially at infinity, the integral defining $b$ is finite and positive, with $b(r)=O(r)$ at zero and $b(r)=O(r^{-5})$ at infinity. In particular, $b\in L^1(0,\infty)$.

To use the polarized identity \eqref{eq3.1} with $w=F_Z$, let $w_{\varepsilon,R}=\chi_{\varepsilon,R}w$, where $\chi_{\varepsilon,R}$ is a smooth cutoff equal to one on $[2\varepsilon,R]$ and supported in $[\varepsilon,2R]$. For $j=1,2$, choose $|\chi_{\varepsilon,R}^{(j)}|\le C_j\varepsilon^{-j}$ on $[\varepsilon,2\varepsilon]$ and $|\chi_{\varepsilon,R}^{(j)}|\le C_jR^{-j}$ on $[R,2R]$. Lemma~\ref{Lem2.3} gives
\[
\frac{w_{\varepsilon,R}'}{z}
=\chi_{\varepsilon,R}Z+\frac r2\chi_{\varepsilon,R}'z
\longrightarrow Z\quad\text{in }\Hh
\quad(\varepsilon\downarrow0,\ R\to\infty).
\]
Indeed, the squared $\Hh$ error near zero is $O(\varepsilon^{\nu-2})$, and the error at infinity is bounded by a polynomial in $R$ times $e^{-2R}$. For fixed $F\in C_c^\infty(0,\infty)$, the local term in the polarized form agrees with its limit once the cutoffs lie outside the support of $F$. The nonlocal pairing converges by the preceding bounds on $B_\nu$, while the left side converges by continuity of $q_\nu$. Since $q_\nu[Z,y_F]=-2\ip{z}{y_F}=0$, we obtain
\begin{equation*}
\mathcal A U=b\qquad\text{on }(0,\infty).
\end{equation*}
This identity holds in distributions. The endpoint estimates give
\[
U(0)=U(\infty)=0,\qquad w(0)=w(\infty)=0.
\]
Indeed, at zero $U=O(r^\nu)$, while at infinity $U=O(e^{-2r})$.

Let
\[
P(r)=\int_0^r\rho(s)\dd s,\qquad P_\infty=\|z\|_2^2,
\qquad
G(r,s)=\frac{P(\min\{r,s\})
(P_\infty-P(\max\{r,s\}))}{P_\infty}.
\]
The functions $1$ and $P$ solve the homogeneous equation for $\mathcal A$. The unit jump of $-\rho^{-1}\partial_r G$ at $r=s$ identifies $G$ as its positive Green kernel with zero endpoint values. Since $G$ is bounded and $b\in L^1$, the first integral below converges and vanishes at both endpoints.

For the second integral, use $U/\rho=O(r)$ at zero, $U/\rho=1+O((1+\log r)/r)$ at infinity and $P_\infty-P(r)=O(e^{-2r})$. Splitting at $s=r$ gives convergence and zero endpoint values, with an $O((1+r)e^{-2r})$ bound at infinity. The integrals solve $\mathcal A U=b$ and $\mathcal A w=U/\rho$. Since no nonzero linear combination of $1$ and $P$ vanishes at both endpoints, it follows that
\begin{equation}\label{eq3.3}
U(r)=\int_0^\infty G(r,s)b(s)\dd s>0,
\qquad
w(r)=\int_0^\infty G(r,s)\frac{U(s)}{\rho(s)}\dd s.
\end{equation}

For $F\in C_c^\infty(0,\infty)$, the same Green representation yields
\[
F(r)=\int_0^\infty G(r,s)(\mathcal A F)(s)\dd s.
\]
Cauchy--Schwarz with the positive weight $U/\rho$ gives
\[
\frac{F(r)^2}{w(r)}
\le\int_0^\infty G(r,s)\frac{\rho(s)}{U(s)}(\mathcal A F(s))^2\dd s.
\]
Multiply by $b(r)$ and integrate. Tonelli's theorem and the first identity in \eqref{eq3.3} imply
\begin{equation*}
\int_0^\infty\frac{b}{w}F^2\dd r
\le\int_0^\infty\rho(\mathcal A F)^2\dd r.
\end{equation*}
By the definition $b(r)=2\int_0^\infty B_\nu(r,s)w(s)\dd s$, subtraction of the nonlocal term in \eqref{eq3.1} and symmetry give \eqref{eq3.2} for $y=y_F$.

For $y\perp z$, take smooth $y_n$ supported in compact annuli and converging to $y$ in $\Hh$. Subtract a multiple of a fixed annular test function to enforce $\ip{z}{y_n}=0$. Then $F_{y_n}$ is compactly supported in $(0,\infty)$ and lies in $H^2$ on its supporting annulus. Approximation there in $H^2$ preserves convergence of $F'/z$ in $H^1$, so the inequality holds for $F_{y_n}$. The form $q_\nu$ is continuous on $\Hh$, and
\[
\sup_{r>0}|F_{y_n}(r)-F_y(r)|
\le\|z\|_2\|y_n-y\|_2.
\]
Thus $F_{y_n}\to F_y$ uniformly, and Fatou's lemma gives \eqref{eq3.2}.

If $\Lcal_\nu y=0$, Lemma~\ref{Lem3.2} gives $y\perp z$ and $q_\nu[y]=0$. The kernel $B_\nu(r,s)w(r)w(s)$ is positive off the diagonal, so equality in \eqref{eq3.2} forces $F_y/w$ to be constant almost everywhere. Local continuity makes it constant everywhere. Thus $y=cZ$, and \eqref{eq2.16} gives $c=0$. Finally, \eqref{eq3.2} makes $q_\nu$ nonnegative on the codimension-one subspace $z^\perp$, so $\ind(\Lcal_\nu)\le1$. Since $q_\nu[z]=-2\int_0^\infty Vz^2\dd r<0$, the index is exactly one.
\end{proof}

We next give the endpoint bounds for a radial kernel element that are needed below.

\begin{lemma}\label{Lem3.4}
Let $z\in\mathcal P_\nu$ and suppose $y\in\Hh$ satisfies $\Lcal_\nu y=0$ and $\ip{z}{y}=0$. With $F_y$ and $w$ as above, set
\[
g=F_y/w,\qquad k=g'.
\]
At zero, $g(r)=g(0)+O(r^2)$, $k(r)=O(r)$ and $k'(r)=O(1)$. At infinity, for some $c\in\R$,
\begin{equation}\label{eq3.4}
\frac{y(r)}{z(r)}=c+O\left(\frac{1+\log r}{r}\right),\qquad
g(r)=O(r^{-1}),\qquad k(r)=O(r^{-1}),\qquad k'(r)=O(r^{-1}).
\end{equation}
\end{lemma}

\begin{proof}
Let $W=\mathcal V_\nu[z,y]$, which is bounded by Lemma~\ref{Lem2.2}. The equation for $r^{-(\nu-1)/2}y$ is regular singular, with bounded forcing $2QW$. The construction in \eqref{eq2.17}, including this forcing, gives a constant plus $O(r^2)$ and a derivative of order $O(r)$. Finite energy excludes the singular homogeneous solution. Substitution into $F_y/w$ gives the origin estimates.

The equation
\[
-y''+\left(1+\frac{\gamma_\nu}{r^2}-V\right)y=2zW
\]
shows that $y''\in L^2(1,\infty)$. Hence $y,y'\to0$ at infinity. Put $v=y/z$ and $\rho=z^2$. The equations for $y$ and $z$ imply $(\rho v')'=-2\rho W$, and
\[
\rho v'=zy'-z'y\longrightarrow0.
\]
By Cauchy--Schwarz in the positive kernel $a_\nu$ and \eqref{eq2.10},
\[
|W(r)|\le\mathcal V_\nu[z](r)^{1/2}
\mathcal V_\nu[y](r)^{1/2}
\le C\frac{1+\log r}{r^2}.
\]
Together with the bound $\int_r^\infty\rho\le C\rho(r)$ from Lemma~\ref{Lem2.3}, this yields
\begin{equation*}
v'(r)=\frac2{\rho(r)}\int_r^\infty\rho(s)W(s)\dd s
=O\left(\frac{1+\log r}{r^2}\right).
\end{equation*}
Thus $v$ has a finite limit $c$ with the stated error. Orthogonality gives
\[
F_y(r)=-\int_r^\infty\rho(s)v(s)\dd s,
\]
so $g=O(r^{-1})$. Finally,
\[
g'=\frac2r(v-\eta g),\qquad
g''=-\frac{g'}r+\frac2r(v'-\eta'g-\eta g').
\]
The bounds for $v,v',\eta,\eta'$ prove the remaining estimates in \eqref{eq3.4}.
\end{proof}

For $3\le\nu<4$, the off-diagonal part of $B_\nu$ is negative. We use this sign only after expressing its distributional action on a squared difference as an ordinary integral.

\begin{lemma}\label{Lem3.5}
Suppose $3\le\nu<4$. Let $g$ be locally $H^2$ on the open half-line and satisfy, for some $N\ge0$,
\begin{equation}\label{eq3.5}
\begin{array}{lll}
g=g_0+O(r^2),& g'=O(r),&g''=O(1)\quad(r\downarrow0),\\
g=O((1+r)^N),&g'=O((1+r)^N),&g''=O((1+r)^N)\quad(r\to\infty).
\end{array}
\end{equation}
The bounds for $g''$ are understood almost everywhere. Put $F=wg$, $y=F'/z$ and $k=g'$. Then $y\in\Hh$, $y\perp z$, and
\begin{equation}\label{eq3.6}
q_\nu[y]=\int_0^\infty\frac{\rho r^2}{4}|k'|^2\dd r
+\int_0^\infty c_z(r)k(r)^2\dd r
-\iint_{(0,\infty)^2}\mathcal K_z(t,u)k(t)k(u)\dd t\dd u,
\end{equation}
where
\begin{equation}\label{eq3.7}
c_z(r)=\rho(r)\left(\frac{\eta(r)+1/2}{2}
-\frac{3r}{2}\eta'(r)\right)
\end{equation}
and, for $t\ne u$,
\begin{equation}\label{eq3.8}
\mathcal K_z(t,u)
=-2\int_0^{\min(t,u)}\int_{\max(t,u)}^\infty
B_\nu(r,s)w(r)w(s)\dd s\dd r>0.
\end{equation}
The double integral in \eqref{eq3.6} is absolutely convergent. Diagonal values of $\mathcal K_z$ may be chosen arbitrarily.
\end{lemma}

\begin{proof}
We first take $g$ smooth with compact support in $(0,\infty)$. The integral formula in Lemma~\ref{Lem3.2} extends the compactly supported pairings to $H^2$ functions on annuli by approximation. Let $w_{\varepsilon,R}=\chi_{\varepsilon,R}w$ be the cutoffs from Proposition~\ref{Prop3.3}. For $G\in H^2$ compactly supported in $(0,\infty)$ we define
\[
\ip{B_\nu}{w\otimes G}
:=\lim_{\varepsilon\downarrow0,\ R\to\infty}
\ip{B_\nu}{w_{\varepsilon,R}\otimes G}.
\]
This limit exists and is independent of the cutoffs: $w_{\varepsilon,R}'/z\to Z$ in $\Hh$, and the nonlocal bilinear form in \eqref{eq2.13} is continuous. In particular, it defines the pairing with $G=wg^2$. Since $F_{y_{wg^2}}=wg^2$ vanishes at both endpoints, \eqref{eq2.16} gives
\[
q_\nu[Z,y_{wg^2}]=0.
\]
Subtract this identity from \eqref{eq3.1} with $F=wg$. The local contribution is
\[
\int_0^\infty\rho\bigl((\mathcal A(wg))^2-(\mathcal A w)\mathcal A(wg^2)\bigr)\dd r.
\]
Using $w'/\rho=\eta$ and $w/\rho=r/2$, its integrand expands as
\[
\frac{\rho r^2}{4}g''^2
+\rho r(\eta+1/2)g'g''
+\rho\bigl((\eta+1/2)^2-r\eta'\bigr)g'^2.
\]
Integrating the middle term by parts gives the local terms in \eqref{eq3.6}, with coefficient \eqref{eq3.7}. Both local terms are unchanged when $k$ is replaced by $|k|$, irrespective of the sign of $c_z$.

By symmetry, the remaining contribution is the extended pairing of $B_\nu$ with
\[
\Phi(r,s)=w(r)w(s)(g(r)-g(s))^2.
\]
More precisely, although $\Phi$ need not be compactly supported, we set
\[
\ip{B_\nu}{\Phi}
:=\lim_{\varepsilon\downarrow0,\ R\to\infty}
\ip{B_\nu}{\chi_{\varepsilon,R}(r)\chi_{\varepsilon,R}(s)\Phi}.
\]
For cutoffs equal to one on the support of $g$, the expression inside the limit is
$2\ip{B_\nu}{w_{\varepsilon,R}\otimes wg^2}
-2\ip{B_\nu}{wg\otimes wg}$, so the limit exists by the preceding definition.
On compact annuli, $\Phi=O(|r-s|^2)$ and its first derivatives are $O(|r-s|)$. There,
\[
|a_\nu|\le C(1+|\log|r-s||),\quad
|\nabla a_\nu|\le C(1+|r-s|^{-1}),\quad
|B_\nu|\le C(1+|r-s|^{-2}).
\]
Integrate by parts outside $|r-s|\le\varepsilon$. On a fixed compact annulus, the boundary terms involving $a_\nu\nabla\Phi$ and $\nabla a_\nu\Phi$ are bounded by $C\varepsilon(1+|\log\varepsilon|)$ and tend to zero. Local integrability also makes the integral of $a_\nu\partial_r\partial_s\Phi$ over the removed strip tend to zero. Hence
\begin{equation}\label{eq3.9}
\ip{B_\nu}{\Phi}
=\iint_{(0,\infty)^2}
B_\nu(r,s)w(r)w(s)(g(r)-g(s))^2\dd r\dd s,
\end{equation}
where the right side uses the ordinary kernel off the diagonal. To remove the annular cutoffs, note that for compactly supported $g$ the integrand vanishes unless at least one radius belongs to its support. Near the other endpoint, $|B_\nu(r,s)|\le C_\nu rs^{-5}$ for $r<s/2$, together with $w(r)=O(r^\nu)$ at zero and exponential decay at infinity, is integrable. Thus dominated convergence identifies the extended pairing with the absolutely convergent integral in \eqref{eq3.9}.

For $r<s$, substitute $g(s)-g(r)=\int_r^s k(t)\dd t$ into \eqref{eq3.9}. Fubini's theorem and symmetry give the last term of \eqref{eq3.6}, with kernel \eqref{eq3.8}. Positivity follows from \eqref{eq2.8}.

To extend the identity, we first prove that \eqref{eq3.5} implies
\begin{equation}\label{eq3.10}
\iint_{r<s}|B_\nu(r,s)|w(r)w(s)
\left(\int_r^s|k(t)|\dd t\right)^2\dd r\dd s<\infty.
\end{equation}
On compact annuli, the squared integral cancels the diagonal singularity. For $0<r<s<1$, combine $w(r)=O(r^\nu)$, homogeneity of $B_\nu$ and $\int_r^s|k|=O(s^2-r^2)$. When $s>2r$, use $|B_\nu(r,s)|\le C_\nu rs^{-5}$; the integral of $|k|$ grows at most polynomially in $s$, while $w(s)$ decays exponentially. For $1<r<s<2r$, use
\[
|B_\nu(r,s)|\le C_\nu r^{-4}
\left(1+\left|\frac{s-r}{r}\right|^{\nu-5}\right),\qquad
\int_r^s|k|\le C(1+s)^N(s-r).
\]
The resulting singularity has order $|s-r|^{\nu-3}$, which is integrable for $\nu\ge3$, and the exponential factor $w(r)w(s)$ controls infinity. This proves \eqref{eq3.10} and justifies Fubini. When $\nu=3$, $\mathcal K_z$ has at most a logarithmic singularity on compact portions of the diagonal.

Finally,
\[
y=z\left(\eta g+\frac r2k\right).
\]
The endpoint estimates for $z,\eta,g,k$ give $y\in\Hh$ and $\ip{z}{y}=F(\infty)-F(0)=0$. Cut off $g$ near zero and infinity, and then mollify on compact subintervals.

Let $\chi_\varepsilon=0$ on $(0,\varepsilon)$ and $\chi_\varepsilon=1$ on $(2\varepsilon,\infty)$, with $|\chi_\varepsilon^{(j)}|\le C_j\varepsilon^{-j}$ for $j=1,2$. The added first and second derivatives of $\chi_\varepsilon g$ are $O(\varepsilon^{-1})$ and $O(\varepsilon^{-2})$. Since $\rho=O(r^{\nu-1})$, the local energy error is $O(\varepsilon^{\nu-2})$, and the corresponding vectors converge in $\Hh$.

For the nonlocal error, put $e_\varepsilon=(\chi_\varepsilon-1)g$ and $H_\varepsilon(r,s)=\int_r^s|e_\varepsilon'(t)|\dd t$ for $r<s$. The origin bounds give $H_\varepsilon\le C\min\{(s-r)/\varepsilon,1\}$, and $H_\varepsilon=0$ if $r\ge2\varepsilon$. On $r<s\le2r$, homogeneity and the diagonal estimate yield
\[
|B_\nu(r,s)|\le C_\nu r^{-4}
\left(1+\left|\frac{s-r}{r}\right|^{\nu-5}\right).
\]
Since $w(r)w(s)\le Cr^{2\nu}$ there,
\[
\iint_{r<s\le2r}|B_\nu|w(r)w(s)H_\varepsilon^2\dd r\dd s
\le C_\nu\varepsilon^{-2}\int_0^{2\varepsilon}r^{2\nu-1}\dd r
\le C_\nu\varepsilon^{2\nu-2}.
\]
On $s>2r$, the separated-radius bound gives
\[
\begin{aligned}
&\iint_{s>2r}|B_\nu|w(r)w(s)H_\varepsilon^2\dd r\dd s\\
&\quad\le C_\nu\int_0^{2\varepsilon}r^{\nu+1}
\int_{2r}^\infty s^{-5}w(s)\dd s\dd r\\
&\quad\le C_\nu\int_0^{2\varepsilon}
r^{\nu+1}(r^{\nu-4}+1)\dd r
\le C_\nu\varepsilon^{2\nu-2},
\end{aligned}
\]
using $w(s)=O(s^\nu)$ at zero and exponential decay at infinity. Both estimates hold for each fixed $3\le\nu<4$. Cauchy--Schwarz with the positive weight $-B_\nu(r,s)w(r)w(s)$ makes the cross terms tend to zero.

For a cutoff $\chi_R$ equal to one below $R$ and zero above $2R$, put $e_R=(\chi_R-1)g$. The derivative $e_R'$ vanishes below $R$ and is bounded by $C(1+r)^N$, uniformly in $R$. Consequently,
\[
\int_r^s|e_R'(t)|\dd t
\le C(1+s)^N(s-r)\boldsymbol 1_{\{s>R\}}.
\]
The bounds just used to prove \eqref{eq3.10} provide an integrable majorant on $s>R>2$, so dominated convergence makes this nonlocal error tend to zero. The local errors and the $\Hh$ error for the associated vectors tend to zero because all cutoff derivatives grow at most polynomially and $z$ decays exponentially. Finally, mollification in $H^2$ on a fixed compact annulus gives uniform convergence of $g$ and $g'$. If $e$ denotes the mollification error, then $\int_r^s|e'|\le\|e'\|_\infty|s-r|$; the same diagonal and endpoint bounds make its nonlocal contribution tend to zero. The local terms converge by the $H^2$ approximation. Passing successively through these three approximations proves \eqref{eq3.6}.
\end{proof}

\begin{proposition}\label{Prop3.6}
Suppose $3\le\nu<4$ and $\ind(\Lcal_\nu)=1$. Then $\ker\Lcal_\nu=\{0\}$.
\end{proposition}

\begin{proof}
Let $y\in\ker\Lcal_\nu$. By Lemma~\ref{Lem3.2}, $y\perp z$ and $y\perp Vz$. Lemma~\ref{Lem3.4} shows that $g=F_y/w$ satisfies \eqref{eq3.5}. Put $k=g'$.

Suppose $k$ has positive and negative parts which are both nonzero. Define
\[
\widetilde g(r)=\int_0^r|k(t)|\dd t,
\qquad
\widetilde y=z\left(\eta\widetilde g+\frac r2|k|\right).
\]
By Lemma~\ref{Lem3.4}, $\widetilde g=O(r^2)$ at zero and $\widetilde g=O(1+\log r)$ at infinity. Since $k\in H^1_{\mathrm{loc}}$, also $|k|\in H^1_{\mathrm{loc}}$ and $(|k|)'=\operatorname{sgn}(k)k'$ almost everywhere. Thus $\widetilde g\in H^2_{\mathrm{loc}}$ satisfies \eqref{eq3.5}, and Lemma~\ref{Lem3.5} gives $\widetilde y\in\Hh$ and $\widetilde y\perp z$. Under $k\mapsto|k|$, the local terms of \eqref{eq3.6} are unchanged because $|(|k|)'|=|k'|$ almost everywhere, while the change in the integral term is
\[
-\iint\mathcal K_z(t,u)
\bigl(|k(t)||k(u)|-k(t)k(u)\bigr)\dd t\dd u<0.
\]
The inequality is strict because the positive and negative sets of $k$ both have positive measure and $\mathcal K_z>0$ away from the diagonal. Consequently,
\[
q_\nu[\widetilde y]<q_\nu[y]=0,
\]
contradicting Lemma~\ref{Lem3.1}.

It follows that $k$ has one sign, so $g$ is monotone. Since $g(\infty)=0$, the functions $g$ and $F_y$ each have one sign. However,
\[
0=\ip{Vz}{y}=\int_0^\infty VF_y'\dd r
=-\int_0^\infty V'F_y\dd r.
\]
The endpoint terms vanish by Lemma~\ref{Lem3.4}. Since $V'<0$, a nonzero function $F_y$ with one sign cannot satisfy this identity. Thus $F_y=0$, and $y=0$.
\end{proof}

Propositions~\ref{Prop3.3} and \ref{Prop3.6}, together with \cite[Theorems~3--4]{KLR2009} at $\nu=4$, give a trivial radial kernel for every $z\in\mathcal P_\nu$ with radial Morse index one. The Fredholm property in Lemma~\ref{Lem2.2} then makes $\Lcal_\nu:\Hh\to\Hh^*$ invertible.

\Needspace{7\baselineskip}
\section{Uniqueness and the full linearized kernel}\label{Sec4}

We first recall the variational characterization in integer dimensions and determine the Morse index. This also applies to the four-dimensional solution from which the continuation starts.

\begin{lemma}\label{Lem4.1}
The least action of a nonzero solution of \eqref{eq1.1} is $j_d$. Every ground state satisfies \eqref{eq1.4} and, after a phase rotation and a translation, is positive, radial and strictly decreasing. At this positive representative, $L_+$ has Morse index one. Every minimizer of $M(u)A(u)/\Dcal(u)$ with $M(u)=A(u)=m_d$ is a ground state of \eqref{eq1.1}.
\end{lemma}

\begin{proof}
For $u\ne0$, the unique positive multiple $tu$ on the Nehari manifold satisfies
\[
t^2=\frac{A(u)+M(u)}{\Dcal(u)},
\qquad S(tu)=\frac{(A(u)+M(u))^2}{4\Dcal(u)}.
\]
The Nehari minimum agrees with the interpolation minimum:
\begin{equation}\label{eq4.1}
\inf_{u\ne0}\frac{(A(u)+M(u))^2}{4\Dcal(u)}
=\inf_{u\ne0}\frac{M(u)A(u)}{\Dcal(u)}=j_d.
\end{equation}
The arithmetic--geometric mean inequality bounds the left-hand infimum below by $j_d$. For the reverse inequality, dilate $u$ so that $\lambda^2A(u)=M(u)$. Attainment follows from \cite{MorozVS2013} and \eqref{eq4.1}.

To verify that the Nehari constraint is natural, set $G(u)=A(u)+M(u)-\Dcal(u)$. On $G(u)=0$, one has $G'(u)u=-2(A(u)+M(u))\ne0$. Testing $S'(u)=\mu G'(u)$ against $u$ forces $\mu=0$, so each Nehari minimizer solves \eqref{eq1.1}. Conversely, every nonzero solution belongs to the Nehari manifold. By testing the equation and differentiating along the mass-preserving dilation, we obtain $A(Q)=M(Q)$ and $\Dcal(Q)=2M(Q)$. Hence $S(Q)=M(Q)/2$, proving \eqref{eq1.4} for ground states.

By \cite[Theorem~3]{MorozVS2013}, a real ground state has one sign and is radial and decreasing about some point. For a complex ground state, $A(|Q|)\le A(Q)$ by the diamagnetic inequality, while $M$ and $\Dcal$ are unchanged. Strict inequality would make the action of the Nehari multiple of $|Q|$ less than $j_d$, contrary to \eqref{eq4.1}. Thus equality holds, and $|Q|$ is a positive Nehari minimizer. Elliptic regularity applied to the real and imaginary parts of $Q$ yields $Q\in C^1(\R^d;\C)$. Since $|Q|>0$ everywhere, the function $\zeta=Q/|Q|$ is globally defined and belongs to $C^1(\R^d;\C)$. The identity
\[
|\nabla Q|^2=|\nabla|Q||^2+|Q|^2|\nabla\zeta|^2
\]
and equality in the diamagnetic inequality imply $\nabla\zeta=0$. Since $\R^d$ is connected, $\zeta$ is constant.

For a positive ground state, \eqref{eq1.5} implies $E(u)\ge0$ on $M(u)=m_d$, with equality at $Q$. Taking the second variation on this mass surface, we obtain
\begin{equation}\label{eq4.2}
\ip{L_+h}{h}\ge0\qquad\text{whenever }\ip{Q}{h}=0.
\end{equation}
The constrained Hessian is $S''(Q)$, since $S=E+M/2$ and $S'(Q)=0$. Since $\ip{L_+Q}{Q}=-2\Dcal(Q)<0$, there is a negative direction. A two-dimensional negative subspace would intersect $Q^\perp$, contradicting \eqref{eq4.2}. Thus the full Morse index is one. Since $Q$ is radial, the radial Morse index is also one.

Finally, a minimizer of the interpolation quotient has the Euler--Lagrange equation
\[
-\Delta u+\frac{A(u)}{M(u)}u
=\frac{2A(u)}{\Dcal(u)}(|x|^{-2}*|u|^2)u.
\]
If $M(u)=A(u)=m_d$, then $\Dcal(u)=2m_d$ by equality in \eqref{eq1.5}. Both coefficients in this equation are one, and its action equals $j_d$.
\end{proof}

For the continuation argument, fix $\nu_*>4$ and set
\[
\mathcal Z_{\nu_*}=\{(\nu,z)\in[3,\nu_*]\times\Hh:
z\in\mathcal P_\nu,\ \ind(\Lcal_\nu)=1\}.
\]
Every profile in $\mathcal P_\nu$ is positive and strictly decreasing by Lemma~\ref{Lem2.3}.

\begin{lemma}\label{Lem4.2}
Every point $(\nu_0,z_0)\in\mathcal Z_{\nu_*}$ has a neighborhood in which the zeros of $\Fcal$ form a unique continuous graph $z=z(\nu)$. This graph remains in $\mathcal Z_{\nu_*}$. At an endpoint of $[3,\nu_*]$, the statement uses a relative, one-sided parameter neighborhood.
\end{lemma}

\begin{proof}
The radial nondegeneracy established in Section~\ref{Sec3}, together with Lemma~\ref{Lem2.2}, shows that
\[
D_z\Fcal(\nu_0,z_0):\Hh\longrightarrow\Hh^*
\]
is an isomorphism. The implicit function theorem with a continuous parameter supplies the unique local graph. At $\nu_0=3$ or $\nu_*$, its contraction proof applies in the relative parameter neighborhood.

We show that the associated profile $Q$ remains positive and strictly decreasing. By Lemma~\ref{Lem2.2},
\begin{equation}\label{eq4.3}
\|\mathcal V_\nu[z(\nu)]-\mathcal V_{\nu_0}[z_0]\|_\infty\longrightarrow0.
\end{equation}
On compact annuli, the radial profiles converge in $C^1$ by the equation and convergence in $\Hh$. Choose $0<r_0<R_0$ so that the reference potential exceeds $1+\varepsilon$ on $[0,r_0]$ and is less than $1/4$ on $[R_0,\infty)$. For sufficiently small $r_0$, the regular Volterra solution with value one at zero stays positive for nearby parameters and potentials. Finite energy excludes the singular solution, as in Lemma~\ref{Lem2.3}. By convergence at $r_0$, $Q(r_0)>0$, so the coefficient of the regular solution is positive. Integrating
\[
(r^{\nu-1}Q')'=(1-V)Qr^{\nu-1}
\]
then shows that $Q'<0$ on $(0,r_0]$.

On $[r_0,R_0]$, positivity and strict decrease follow from $C^1$ convergence. On $[R_0,\infty)$, we have $1-V>1/2$ by \eqref{eq4.3}. The negative part $z_- =\max\{-z,0\}$ has zero trace at $R_0$. Testing the conjugated equation against $-z_-$ yields
\[
\int_{R_0}^\infty\left(|z_-'|^2+
\left[1+\frac{\gamma_\nu}{r^2}-V\right]z_-^2\right)\dd r=0.
\]
Since $\gamma_\nu\ge0$, this forces $z_-=0$. The strong maximum principle ensures $Q>0$. Exterior comparison implies exponential decay, so we may integrate from infinity to obtain
\[
r^{\nu-1}Q'(r)
=-\int_r^\infty (1-V(s))Q(s)s^{\nu-1}\dd s<0.
\]
After conjugation, the test functions used on $(R_0,\infty)$ belong to $H_0^1(R_0,\infty)$; the energy bounds therefore justify these integrations.

Represent the linearized forms by bounded self-adjoint operators on $\Hh$ using its Riesz map. These operators vary continuously in norm and are compact perturbations of uniformly positive operators. Since the operator at $(\nu_0,z_0)$ is invertible, its spectrum has a gap around zero. The negative spectral projection therefore retains its finite rank for nearby parameters, so the Morse index remains one.
\end{proof}

\begin{proposition}\label{Prop4.3}
For every real $\nu\ge3$, equation \eqref{eq2.2} has exactly one solution $z$ such that $Q=r^{-(\nu-1)/2}z$ is positive and strictly decreasing and $\ind(\Lcal_\nu)=1$. For every real $\nu\ge4$, the set $\mathcal P_\nu$ consists of exactly one element, without an index assumption.
\end{proposition}

\begin{proof}
Fix $\nu_*>4$. By Proposition~\ref{Prop2.4} and Lemma~\ref{Lem2.5}, every sequence in $\mathcal Z_{\nu_*}$ has a strongly convergent subsequence whose limit is a nonzero solution. Local uniform convergence of the radial profiles preserves nonnegativity and monotonicity. The limit is strictly positive and strictly decreasing by Lemma~\ref{Lem2.3}.

Operator-norm convergence implies that the limiting Morse index is at most one. Indeed, a two-dimensional negative subspace for the limiting form would remain negative for the approximating forms when $n$ is large. On the other hand,
\[
q_\nu[z]=-2\int_0^\infty Vz^2\dd r<0,
\]
so the limiting index is one. The radial kernel at this limit is trivial by Section~\ref{Sec3}. Thus $\mathcal Z_{\nu_*}$ is compact, and Lemma~\ref{Lem4.2} applies at every point.

Let $N(\nu)$ be the number of solutions at parameter $\nu$ in $\mathcal Z_{\nu_*}$. For fixed $\nu$, these solutions form a compact set of isolated points by Lemma~\ref{Lem4.2}, so $N(\nu)$ is finite. Fix $\nu_0$ and apply that lemma at each solution. Any additional solutions at parameters tending to $\nu_0$ would, by compactness, converge along a subsequence to a solution at $\nu_0$, contrary to local uniqueness. The same argument shows that if $N(\nu_0)=0$, then $N(\nu)=0$ nearby. Hence $N$ is locally constant on $[3,\nu_*]$.

At $\nu=4$, a positive radial ground state exists by Lemma~\ref{Lem4.1}. The uniqueness theorem \cite[Theorem~3]{KLR2009} applies to all positive radial $H^1$ solutions, so every positive radial solution is this ground state. Its radial Morse index is one by Lemma~\ref{Lem4.1}, and hence $N(4)=1$. Since $[3,\nu_*]$ is connected, $N(\nu)=1$ throughout this interval. The choice of $\nu_*>4$ was arbitrary, so the assertion holds for every $\nu\ge3$. For $\nu>4$, Proposition~\ref{Prop3.3} shows that every $z\in\mathcal P_\nu$ has radial Morse index one. At $\nu=4$, the same conclusion follows from the uniqueness theorem just cited. Thus $\mathcal P_\nu$ is a singleton for every $\nu\ge4$.
\end{proof}

We next pass from radial nondegeneracy to the full kernel. Fix a positive ground state $Q$ in $\R^d$, and let $q_+$ and $q_-$ denote the forms of $L_+$ and $L_-$ on real $H^1(\R^d)$.

\begin{lemma}\label{Lem4.4}
The operators $L_\pm$ are self-adjoint on $H^2(\R^d)$ and have essential spectrum $[1,\infty)$. Their forms are compact perturbations of the $H^1$ inner product, and every weak $H^1$ kernel element belongs to $H^2$.
\end{lemma}

\begin{proof}
The potential $V_Q$ is bounded and tends to zero. The exchange operator
\[
Kh(x)=2Q(x)\int_{\R^d}|x-y|^{-2}Q(y)h(y)\dd y
\]
is compact on $L^2$. Indeed, first remove $|x-y|<\varepsilon$. Since $Q$ is bounded, the Schur bound for this part is $O(\varepsilon^{d-2})$. Next truncate both variables to a ball. Factoring out half the exponential decay in the variable outside the ball, and splitting the remaining convolution into $|x-y|<1$ and $|x-y|\ge1$, makes the Schur bound for the discarded tails tend to zero. The remaining bounded kernel has compact support and defines a Hilbert--Schmidt operator. Thus $K$ is a norm limit of compact operators.

Truncation and the Rellich theorem show that multiplication by $V_Q$ is compact from $H^1$ to $L^2$. The asserted self-adjointness, essential spectrum and compactness of the perturbations of the forms follow. A weak $H^1$ kernel element has Laplacian in $L^2$, and therefore belongs to $H^2$.
\end{proof}

For an $L^2(\mathbb S^{d-1})$-normalized spherical harmonic $Y_{\ell,m}$ of degree $\ell$, define $a_{d,\ell}$ by
\[
\int_{\mathbb S^{d-1}}
\frac{Y_{\ell,m}(\theta)}{|r\omega-s\theta|^2}\dd\theta
=a_{d,\ell}(r,s)Y_{\ell,m}(\omega).
\]
These coefficients are independent of $m$. In particular, $a_{d,0}=a_d$, and $a_{d,1}$ is the first angular moment in \eqref{eq2.1} at $\nu=d$.

\begin{lemma}\label{Lem4.5}
For $r,s>0$ with $r\ne s$, the angular kernels satisfy
\[
a_{d,\ell}(r,s)>0,\qquad a_{d,\ell+1}(r,s)<a_{d,\ell}(r,s),\qquad \ell\ge0.
\]
\end{lemma}

\begin{proof}
Using $|x-y|^{-2}=\int_0^\infty e^{-t|x-y|^2}\dd t$ in the Funk--Hecke formula \cite{AtkinsonHan2012}, we obtain
\begin{equation}\label{eq4.4}
a_{d,\ell}(r,s)=(2\pi)^{d/2}
\int_0^\infty e^{-t(r^2+s^2)}
(2trs)^{1-d/2}I_{\ell+d/2-1}(2trs)\dd t,
\end{equation}
where $I_\alpha$ is the modified Bessel function. The integral converges for $r\ne s$, and its integrand is positive.

To compare consecutive orders, use the recurrence relations for $I_\alpha$. For $R_\alpha(x)=I_{\alpha+1}(x)/I_\alpha(x)$,
\[
R_\alpha'=1-\frac{2\alpha+1}{x}R_\alpha-R_\alpha^2,
\qquad R_\alpha(x)\sim\frac{x}{2\alpha+2}\quad(x\downarrow0).
\]
For $\alpha>-1/2$, the positive function $R_\alpha$ starts below one. If $x_0>0$ were the first point at which $R_\alpha(x_0)=1$, then $R_\alpha'(x_0)\ge0$, whereas the Riccati equation would imply
\[
R_\alpha'(x_0)=-\frac{2\alpha+1}{x_0}<0.
\]
Consequently $I_{\alpha+1}(x)<I_\alpha(x)$ for $x>0$. In \eqref{eq4.4}, the relevant orders satisfy $\alpha=\ell+d/2-1\ge1/2$, so the strict kernel comparison follows.
\end{proof}

Let $\lambda_{d,\ell}=\ell(\ell+d-2)$. In degree $\ell$ the radial form is
\begin{equation}\label{eq4.5}
\begin{split}
q_\ell[f]={}&\int_0^\infty
\left(f'^2+\left[1-V_Q+\frac{\lambda_{d,\ell}}{r^2}\right]f^2\right)r^{d-1}\dd r\\
&-2\iint a_{d,\ell}(r,s)Q(r)Q(s)f(r)f(s)
r^{d-1}s^{d-1}\dd r\dd s.
\end{split}
\end{equation}
Its form domain consists of the radial coefficients of $H^1$ functions in the spherical harmonic subspace of degree $\ell$.

\begin{lemma}\label{Lem4.6}
The degree-one kernel is generated by $Q'$, and the kernels in every degree $\ell\ge2$ are trivial.
\end{lemma}

\begin{proof}
By translation invariance, $L_+(\partial_{x_j}Q)=0$. Hence $\psi=-Q'>0$ solves the degree-one radial equation. For $f$ smooth and supported in an annulus, set $g=f/\psi$ and test that equation against $\psi g^2$. Subtracting the resulting identity from \eqref{eq4.5} and symmetrizing, we find
\begin{equation*}
\begin{split}
q_1[f]={}&\int_0^\infty\psi^2|g'|^2r^{d-1}\dd r\\
&+\iint a_{d,1}(r,s)Q(r)Q(s)\psi(r)\psi(s)
|g(r)-g(s)|^2r^{d-1}s^{d-1}\dd r\dd s.
\end{split}
\end{equation*}
Both terms are nonnegative. Annular approximation and Fatou's lemma extend the lower bound to the full degree-one form domain. At zero, the cutoff error is bounded by $C\int_0^{2\varepsilon}f^2r^{d-3}\dd r$, which tends to zero by the form-domain condition. At infinity, the corresponding error is bounded by $CR^{-2}\int_R^{2R}f^2r^{d-1}\dd r$ and also tends to zero. On compact annuli, $\psi$ is bounded away from zero, so the quotients converge in $H^1$. For a kernel element, the form vanishes, forcing $f/\psi$ to be constant.

For $\ell\ge2$, replacing $f$ by $|f|$ does not increase the local energy or the negative integral term, since $a_{d,\ell}>0$. It follows from Lemma~\ref{Lem4.5} that
\begin{equation*}
q_\ell[f]\ge q_\ell[|f|]
\ge q_1[|f|]+(\lambda_{d,\ell}-\lambda_{d,1})
\int_0^\infty f^2r^{d-3}\dd r.
\end{equation*}
The right side is strictly positive for nonzero $f$, so the kernel is trivial in every degree $\ell\ge2$.
\end{proof}

\begin{proof}[Proof of Theorem~\ref{Thm1.1}]
By Lemma~\ref{Lem4.1}, every positive ground state centered at the origin corresponds to a half-line solution in the class of Proposition~\ref{Prop4.3}. That proposition proves uniqueness. Lemma~\ref{Lem4.1} also supplies the Morse index and the phase and translation classification of complex ground states.

The radial kernel is trivial by Section~\ref{Sec3}. Spherical harmonic projections are bounded on $H^1$ and commute with $L_+$ by radial symmetry of the potential and interaction. Each projected kernel element therefore belongs to the corresponding radial form domain and satisfies its weak equation. Lemma~\ref{Lem4.6} identifies the remaining kernel elements as the degree-one translation derivatives. This proves the full kernel characterization for $L_+$.

Finally, $L_-Q=0$ and $Q>0$. The ground-state identity
\[
q_-[h]=\int_{\R^d}Q^2|\nabla(h/Q)|^2
\]
initially holds for compactly supported smooth $h$. By approximation, the corresponding lower bound holds for all $h\in H^1$. A zero mode must have constant quotient $h/Q$, and hence is a multiple of $Q$. The converse follows from the equation.

To treat all positive solutions when $d\ge4$, apply the symmetry theorem of Ma and Zhao \cite[Theorem~2]{MaZhao2010}, in the formulation of \cite[Assumption~1.2 and Theorem~1.4]{Li2022}. That formulation uses $\alpha$ for the exponent of $|x|^{-\alpha}$. Sobolev interpolation gives $Q\in L^{2d/(d-1)}$, and the base parameter range holds because $(d-2)/(2d-2)<1/2<d/(2d-2)$. For $\alpha=p=2$, its auxiliary exponents may be chosen as
\[
r=r_1=r_2=r_3=\frac{2d}{d-1},\qquad
t=t_1=d,\qquad s=\frac{2d}{d+1}.
\]
Indeed, $2<r<2d/(d-2)$, $0<1/t=1/d<2/d$, and
\[
\frac1s=\frac1t+\frac1r,\qquad
\frac1t+\frac{d-2}{d}=\frac2r,\qquad
\frac2d\le\frac1s\le\min\left\{1,\frac1r+\frac2d\right\}.
\]
These verify the exponent conditions for every $d\ge3$. Hence every positive $H^1$ solution is radial and decreasing about some point. After translation and conjugation it belongs to $\mathcal P_d$. For $d\ge4$, Proposition~\ref{Prop4.3} identifies it with $Q_d$.
\end{proof}

\section{Quantitative stability and threshold dynamics}\label{Sec5}

We first introduce the correlation with the ground-state orbit. Put
\[
W=(1-\Delta)Q_d=(|x|^{-2}*Q_d^2)Q_d,
\qquad
\mathcal C(u)=\sup_{y\in\R^d}
\left|\ip{u}{Q_d(\cdot-y)}_{H^1}\right|.
\]
The distance to the orbit has an exact formula in terms of $\mathcal C$. We use it to choose a nearest ground state and prove coercivity under the mass, gradient and symmetry constraints. Compactness of normalized extremizing sequences then extends the local estimate to the full constraint set. 

\begin{lemma}\label{Lem5.1}
For $u\in\mathcal N$,
\begin{equation}\label{eq5.1}
\dist_{H^1}(u,\mathcal O_d)^2=4m_d-2\mathcal C(u).
\end{equation}
For every nonzero $u\in H^1(\R^d;\C)$, the supremum defining $\mathcal C(u)$ is positive and attained. If $1<p<\infty$ and $u,v\in H^1\cap L^p$, then
\begin{equation}\label{eq5.2}
|\mathcal C(u)-\mathcal C(v)|
\le\|u-v\|_p\|W\|_{p'},
\qquad p'=\frac p{p-1}.
\end{equation}
In particular, every $u\in\mathcal N$ has a nearest element of $\mathcal O_d$.
\end{lemma}

\begin{proof}
The regularity and decay of the ground state imply that $W$ is bounded and exponentially decaying, so $W\in L^q$ for $1\le q\le\infty$. The equation and local elliptic estimates also give exponential decay of the finitely many derivatives of $Q_d$ and $W$ used below; see \cite{MorozVS2013}. Integration by parts yields
\[
\ip{u}{Q_d(\cdot-y)}_{H^1}
=\int_{\R^d}u(x)W(x-y)\dd x.
\]
This function of $y$ is continuous and tends to zero at infinity. Indeed, approximate both factors in $L^2$ by compactly supported smooth functions; Cauchy--Schwarz bounds the approximation error uniformly in $y$.

The correlation cannot vanish identically for nonzero $u$. Exponential decay of $W$ makes its Fourier transform real analytic on $\R^d$, and $\widehat W(0)=\int W>0$. The zero set of this nontrivial real analytic function has Lebesgue measure zero. Taking Fourier transforms of a vanishing correlation would therefore imply $\widehat u=0$ almost everywhere. Its absolute value consequently has a positive maximum. Expanding the squared $H^1$ norm and minimizing first in phase proves \eqref{eq5.1}, since $\|u\|_{H^1}^2=\|Q_d\|_{H^1}^2=2m_d$. H\"older's inequality proves \eqref{eq5.2}.
\end{proof}

\begin{lemma}\label{Lem5.2}
There is $c>0$ such that
\begin{equation}\label{eq5.3}
q_+[a]\ge c\|a\|_{H^1}^2
\end{equation}
for every real $a\in H^1$ satisfying
\begin{equation}\label{eq5.4}
\ip{Q_d}{a}=0,\qquad
\ip{\nabla Q_d}{\nabla a}=0,\qquad
\ip{a}{\partial_{x_j}Q_d}_{H^1}=0\quad(1\le j\le d).
\end{equation}
Also,
\begin{equation}\label{eq5.5}
q_-[b]\ge c\|b\|_{H^1}^2
\qquad\text{if }\ip{b}{Q_d}_{H^1}=0.
\end{equation}
The same coercivity conclusions hold with $L^2$ orthogonality to the symmetry directions.
\end{lemma}

\begin{proof}
By \eqref{eq4.2}, $q_+$ is nonnegative on $Q_d^\perp$, where this perpendicular is taken in $L^2$. If $a\perp Q_d$ and $q_+[a]=0$, polarization on this subspace implies
\[
L_+a=c_0Q_d\quad\text{in }H^{-1}
\]
for some $c_0\in\R$. Since $L_+\Lambda Q_d=-2Q_d$, Theorem~\ref{Thm1.1} gives
\[
a=-\frac{c_0}{2}\Lambda Q_d+\sum_{j=1}^dc_j\partial_{x_j}Q_d.
\]
Differentiating $A(\lambda^{d/2}Q_d(\lambda\cdot))=\lambda^2m_d$, we obtain
\[
\ip{\nabla Q_d}{\nabla\Lambda Q_d}=m_d.
\]
The gradient constraint in \eqref{eq5.4} forces $c_0=0$. The $H^1$ Gram matrix of the translation derivatives is positive definite, so the remaining conditions force every $c_j$ to vanish.

If \eqref{eq5.3} failed, there would be $a_n$ satisfying \eqref{eq5.4} with $\|a_n\|_{H^1}=1$ and $q_+[a_n]\to0$. Pass to a weak limit $a$. The constraints persist, and the compact terms in the form converge by Lemma~\ref{Lem4.4}. Hence
\[
0\le q_+[a]\le\liminf_nq_+[a_n]=0,
\]
so $a=0$ by the preceding null-space calculation. Compactness would then imply $q_+[a_n]=1+o(1)$, a contradiction. The proof of \eqref{eq5.5} is the same, using $q_-\ge0$ and $\ker L_-=\spanop\{Q_d\}$.

For the last assertion, project an $L^2$-orthogonal perturbation onto the corresponding $H^1$ complement by subtracting a linear combination of the symmetry modes. The quadratic form is unchanged because these modes lie in its kernel. The original $L^2$ orthogonality bounds the coefficients by the $H^1$ norm of the projected vector, through the positive $L^2$ Gram matrix. Thus the two complement norms are equivalent. The translation projection preserves both real normalization constraints, since their pairings with each $\partial_{x_j}Q_d$ vanish.
\end{proof}

\begin{lemma}\label{Lem5.3}
If $M(u_n)=A(u_n)=m_d$ and $\Dcal(u_n)\to2m_d$, then
\[
\dist_{H^1}(u_n,\mathcal O_d)\longrightarrow0.
\]
\end{lemma}

\begin{proof}
Since $2<2d/(d-1)<2d/(d-2)$, estimate \eqref{eq1.3} bounds $\|u_n\|_{2d/(d-1)}$ away from zero. By the vanishing lemma from concentration-compactness theory \cite{Lions1984}, there are a subsequence and centers $y_n$ such that $\int_{B_1(y_n)}|u_n|^2\ge\varepsilon>0$. After translation by $y_n$ and further extraction, local Rellich compactness implies
\[
u_n\rightharpoonup u\ne0\quad\text{in }H^1,
\qquad u_n\to u\quad\text{locally in }L^{2d/(d-1)}
\quad\text{and almost everywhere}.
\]
Set $v_n=u_n-u$. Expanding the squared norms and using weak convergence, we obtain
\begin{equation}\label{eq5.6}
M(u_n)=M(u)+M(v_n)+o(1),\qquad
A(u_n)=A(u)+A(v_n)+o(1).
\end{equation}
We also have the nonlocal splitting
\begin{equation}\label{eq5.7}
\Dcal(u_n)=\Dcal(u)+\Dcal(v_n)+o(1).
\end{equation}
Local $L^{2d/(d-1)}$ convergence, together with approximation of $u$ by compactly supported functions in that space, implies $u\overline{v_n}\to0$ in $L^{d/(d-1)}$. Also, $|v_n|^2$ is bounded in $L^{d/(d-1)}$ and tends to zero almost everywhere, hence converges weakly to zero in $L^{d/(d-1)}$. The Hardy--Littlewood--Sobolev inequality makes the bilinear form
\[
\mathcal B(f,g)=\iint |x-y|^{-2}f(x)g(y)\dd x\dd y
\]
continuous on $L^{d/(d-1)}\times L^{d/(d-1)}$. Expand $|u_n|^2=|u|^2+|v_n|^2+2\Rea(u\overline{v_n})$. Strong convergence eliminates the terms containing the last summand, while $\mathcal B(|u|^2,|v_n|^2)\to0$ by weak convergence. This proves \eqref{eq5.7}.

We have $M(u)>0$ and $A(u)>0$, since a nonzero $L^2$ function cannot have zero gradient. Combining \eqref{eq5.6}--\eqref{eq5.7} with the sharp inequality, we find
\[
2m_d\le\frac{2}{m_d}\bigl[M(u)A(u)+(m_d-M(u))(m_d-A(u))\bigr]\le2m_d.
\]
The difference between $m_d^2$ and the expression in brackets is
\[
M(u)(m_d-A(u))+A(u)(m_d-M(u)).
\]
Both terms are nonnegative. Since $M(u),A(u)>0$, equality forces $M(u)=A(u)=m_d$. The convergence is therefore strong in $H^1$, and the limit is a normalized extremizer. Lemma~\ref{Lem4.1} and Theorem~\ref{Thm1.1} place this limit in $\mathcal O_d$. If the distance did not tend to zero along the full sequence, the same argument applied to a subsequence on which it is bounded away from zero would give a contradiction.
\end{proof}

\begin{proof}[Proof of Theorem~\ref{Thm1.4}]
By Lemma~\ref{Lem5.1}, choose a nearest element $e^{i\theta}Q_d(\cdot-y)$ and write
\[
h=e^{-i\theta}u(\cdot+y)-Q_d=a+ib.
\]
The case $h=0$ is immediate. At a nearest point, differentiation in phase and translation gives
\[
\ip{b}{Q_d}_{H^1}=0,
\qquad
\ip{a}{\partial_{x_j}Q_d}_{H^1}=0,
\qquad
\|h\|_{H^1}=\dist_{H^1}(u,\mathcal O_d).
\]
We first treat small $\|h\|_{H^1}$.

The exact constraints $M(u)=A(u)=m_d$ imply
\begin{equation*}
\ip{Q_d}{a}=-\frac12\|h\|_2^2,
\qquad
\ip{\nabla Q_d}{\nabla a}=-\frac12\|\nabla h\|_2^2.
\end{equation*}
Subtract a linear combination of $Q_d$ and $\Lambda Q_d$ from $a$ to obtain $\widehat a$ satisfying the first two conditions in \eqref{eq5.4}. The constraint matrix on $Q_d,\Lambda Q_d$ is
\[
\begin{pmatrix}m_d&0\\m_d&m_d\end{pmatrix}.
\]
It is invertible, so the coefficients of this correction are $O(\|h\|_{H^1}^2)$. Radial symmetry of $Q_d$ and $\Lambda Q_d$ preserves translation orthogonality. By Lemma~\ref{Lem5.2},
\begin{equation}\label{eq5.8}
q_+[a]+q_-[b]
\ge c\|h\|_{H^1}^2-C\|h\|_{H^1}^3.
\end{equation}

By the bilinear Hardy--Littlewood--Sobolev estimate, the interaction is a continuous quartic polynomial on complex $H^1$ viewed as a real Hilbert space. Since $S'(Q_d)=0$, Taylor expansion yields
\[
S(Q_d+h)-S(Q_d)
=\frac12\bigl(q_+[a]+q_-[b]\bigr)
+O(\|h\|_{H^1}^3+\|h\|_{H^1}^4).
\]
Since the masses agree and $E(Q_d)=0$,
\[
2m_d-\Dcal(u)=4E(u)=4(S(u)-S(Q_d)).
\]
For sufficiently small $\|h\|_{H^1}$, \eqref{eq5.8} proves \eqref{eq1.6} locally.

Fix the neighborhood where the local estimate holds. By Lemma~\ref{Lem5.3}, the deficit has a positive lower bound on the normalized constraint set outside this neighborhood. Otherwise a sequence with deficit tending to zero would have distance to $\mathcal O_d$ tending to zero. On the whole constraint set, Lemma~\ref{Lem5.1} gives
\[
\dist_{H^1}(u,\mathcal O_d)^2\le4m_d.
\]
Combining this bound with the local estimate and the positive lower bound away from $\mathcal O_d$ proves \eqref{eq1.6}.
\end{proof}

\begin{proof}[Proof of Corollary~\ref{Cor1.5}]
Since $u\in H^1$ and $M(u)=m_d>0$, its gradient norm is nonzero and $\lambda$ is well defined. The scaling identities imply
\[
M(v)=A(v)=m_d,
\qquad E(v)=\lambda^2E(u),
\qquad \lambda^2=\frac{m_d}{A(u)}.
\]
By Theorem~\ref{Thm1.4},
\[
\dist_{H^1}(v,\mathcal O_d)^2
\le\frac{1}{c_d}\bigl(2m_d-\Dcal(v)\bigr)
=\frac{4m_d}{c_d}\frac{E(u)}{A(u)},
\]
which is \eqref{eq1.7}. The sharp inequality ensures $E(u)\ge0$ when $M(u)=m_d$. Thus bounded energy and $A(u_n)\to\infty$ imply convergence to $\mathcal O_d$ after normalization.

For the final assertion, define, for all sufficiently large $n$,
\[
c_n=\left(\frac{m_d}{M(u_n)}\right)^{1/2},
\qquad \lambda_n=\left(\frac{M(u_n)}{A(u_n)}\right)^{1/2},
\qquad v_n=c_n\lambda_n^{d/2}u_n(\lambda_n\cdot).
\]
Then $M(v_n)=A(v_n)=m_d$ and
\[
2m_d-\Dcal(v_n)
=2m_d-\frac{m_d^2\Dcal(u_n)}{M(u_n)A(u_n)}
\longrightarrow0,
\]
since $\Dcal(u_n)=2A(u_n)-4E(u_n)\to2m_d$. Applying Theorem~\ref{Thm1.4} directly, we obtain
\[
\dist_{H^1}(v_n,\mathcal O_d)^2
\le\frac{1}{c_d}\left(
2m_d-\frac{m_d^2\Dcal(u_n)}{M(u_n)A(u_n)}\right)
\longrightarrow0.
\]
Choose phases $\theta_n$ and centers $y_n$ for which
\[
\|v_n-e^{i\theta_n}Q_d(\cdot-y_n)\|_{H^1}\longrightarrow0.
\]
The inverse normalization $T_n h=c_n^{-1}\lambda_n^{-d/2}h(\cdot/\lambda_n)$ satisfies
\[
\|T_nh\|_{H^1}^2
=\frac{M(u_n)}{m_d}\|h\|_2^2
+\frac{A(u_n)}{m_d}\|\nabla h\|_2^2.
\]
These operators are therefore uniformly bounded, and
\[
\left\|u_n-c_n^{-1}\lambda_n^{-d/2}e^{i\theta_n}
Q_d(\cdot/\lambda_n-y_n)\right\|_{H^1}\longrightarrow0.
\]
The comparison function has center $\lambda_n y_n$. Since $c_n,\lambda_n\to1$, continuity of dilation at the fixed profile $Q_d$, together with translation invariance of the norm, implies that its difference from $e^{i\theta_n}Q_d(\cdot-\lambda_n y_n)$ tends to zero in $H^1$. This proves the final assertion.
\end{proof}

We now turn to the threshold dynamics.

\begin{proof}[Proof of Corollary~\ref{Cor1.3}]
Li and Zhang \cite[Theorem~1.4]{LiZhang2009} and Gao and Wang \cite{GaoWang2017} assume uniqueness of the ground state in dimensions $d\ge5$ and $d=3$, respectively. Their proofs use this assumption to identify the optimizers of the sharp interpolation inequality. If $M(f)=m_d$ and $E(f)=0$, then $\Dcal(f)=2A(f)$ and
\[
\frac{M(f)A(f)}{\Dcal(f)}=\frac{m_d}{2}=j_d.
\]
After a mass-preserving dilation making the kinetic energy equal to $m_d$, Lemma~\ref{Lem4.1} identifies the rescaled function as a ground state. Theorem~\ref{Thm1.1} therefore gives the phase, translation and dilation classification needed in these variational steps.

Gao and Wang also assume the kernel property in \cite{GaoWang2017}. To verify the spectral condition in their Lemma~1.5, let $\omega$ be an $L^2$-normalized eigenfunction for the unique negative eigenvalue $e_0$ of $L_+$. Lemma~\ref{Lem4.4} and Theorem~\ref{Thm1.1} show that zero is isolated from the positive spectrum and that its eigenspace is spanned by the translation derivatives. The spectral theorem gives, for some $\sigma>0$,
\[
q_+[h]\ge\sigma\|h\|_2^2
\quad\text{if }\ip{h}{\omega}=\ip{h}{\partial_{x_j}Q_d}=0
\quad(1\le j\le d).
\]
Together with the $L^2$-orthogonality version of the estimate for $q_-$ in Lemma~\ref{Lem5.2}, this gives both bounds in \cite[Lemma~1.5]{GaoWang2017}. The $L^2$-bounded perturbations in Lemma~\ref{Lem4.4} also imply $\|h\|_{H^1}^2\le q_+[h]+C\|h\|_2^2$, which upgrades the displayed estimate to an $H^1$ bound. The negative eigenfunction is radial: the radial subspace contains the negative direction $Q_d$, and the full negative eigenspace is one-dimensional. Moreover, \eqref{eq4.2} implies $\ip{Q_d}{\omega}\ne0$, so
\[
e_0\ip{\Lambda Q_d}{\omega}=-2\ip{Q_d}{\omega}\ne0.
\]
This verifies the nonzero pairing used in their modulation argument.

For $d\ge5$, Li and Zhang \cite[Theorem~1.4]{LiZhang2009} require
\[
\|u\|_{L_{t,x}^{2(d+2)/d}(I\times\R^d)}=\infty,
\]
where $I$ is the maximal lifespan; the discussion following their theorem allows this divergence to occur in either time direction. For $d=3,4$, Gao and Wang \cite[Theorem~1.3]{GaoWang2017} classify global solutions that fail to scatter in at least one time direction. Failure of $L^2$ scattering in either time direction forces the corresponding half-line spacetime norm to be infinite by the local scattering criterion. Since $I=\R$, the hypothesis of Li and Zhang follows as well. Thus both classification results apply. Global existence on $\R$ excludes the finite-time pseudoconformal alternative in the former theorem. The remaining solution is a radial solitary wave, whose frequency is $\lambda^2$ under the mass-critical scaling.
\end{proof}

\begin{proof}[Proof of Theorem~\ref{Thm1.2}]
We adapt the classification argument of \cite{Merle1993,HmidiKeraani2005,MXZ2010}, using Corollary~\ref{Cor1.5} for the compactness step. Put $E_0=E(u(0))\ge0$. For a real $C^1$ function $\psi$ with bounded gradient, the sharp inequality gives $E(e^{is\psi}u(t))\ge0$ for every $s\in\R$. Expanding this quadratic polynomial in $s$ and taking its discriminant, as in \cite[Lemma~2.1]{Banica2004} and \cite[Lemma~3.3]{MXZ2010}, yields
\begin{equation}\label{eq5.9}
\left|\int_{\R^d}\nabla\psi\cdot\operatorname{Im}(\overline u\nabla u)\right|^2
\le 2E_0\int_{\R^d}|\nabla\psi|^2|u|^2.
\end{equation}

Choose $t_n\uparrow T$. The $H^1$ local theory and blow-up alternative follow from the admissible pair $(4,2d/(d-1))$ and the Hardy--Littlewood--Sobolev estimate, in every $d\ge3$. In particular, $A(u(t_n))\to\infty$. By Corollary~\ref{Cor1.5}, there exist $\rho_n\to0$, $x_n\in\R^d$ and $\theta_n\in\R$ such that
\[
e^{-i\theta_n}\rho_n^{d/2}u(t_n,x_n+\rho_n x)\longrightarrow Q_d(x)
\quad\text{in }H^1(\R^d).
\]
In particular, for every bounded uniformly continuous function $\psi$,
\begin{equation}\label{eq5.10}
\int_{\R^d}\psi(x)|u(t_n,x)|^2\dd x-m_d\psi(x_n)\longrightarrow0.
\end{equation}
Indeed, strong $L^2$ convergence gives strong $L^1$ convergence of the squared moduli, and uniform continuity applies after restricting the rescaled variable to a fixed ball; the $L^2$ tail of $Q_d$ controls its complement.

After passing to a subsequence, either $x_n\to x_*$ or $|x_n|\to\infty$. In the first case set $a=x_*$, and in the second set $a=0$. Take a nonnegative $\phi\in C_c^\infty(\R^d)$ with $\phi(x)=|x|^2$ for $|x|\le1$ and $|\nabla\phi|^2\le C\phi$, and define
\[
\phi_R(x)=R^2\phi\bigl((x-a)/R\bigr),
\qquad g_R(t)=\int_{\R^d}\phi_R(x)|u(t,x)|^2\dd x.
\]
The continuity equation and \eqref{eq5.9} give
\[
|g_R'(t)|
=2\left|\int\nabla\phi_R\cdot\operatorname{Im}(\overline u\nabla u)\right|
\le C E_0^{1/2}g_R(t)^{1/2},
\]
with $C$ independent of $R$. By \eqref{eq5.10}, $g_R(t_n)\to0$ in either case. Integrating the inequality for $(g_R+\varepsilon)^{1/2}$ and then letting $n\to\infty$ and $\varepsilon\downarrow0$, we obtain $g_R(t)\le CE_0(T-t)^2$. Fatou's lemma therefore gives
\begin{equation}\label{eq5.11}
V_a(t):=\int_{\R^d}|x-a|^2|u(t,x)|^2\dd x
\le CE_0(T-t)^2.
\end{equation}
The alternative $|x_n|\to\infty$ is impossible, since the rescaled convergence gives
\[
\int_{|x-x_n|<1}|u(t_n,x)|^2\dd x\longrightarrow m_d,
\]
whereas \eqref{eq5.11} with $a=0$ would force this integral to tend to zero. Hence $x_n\to x_*$, and \eqref{eq5.11} holds with $a=x_*$.

We next justify the full variance identity at this regularity. For $K(x)=|x|^{-2}$, the localized virial identity reads
\[
\begin{aligned}
g_R''(t)={}&4\operatorname{Re}\int_{\R^d}
D^2\phi_R[\nabla u,\nabla\overline u]
-\int_{\R^d}\Delta^2\phi_R|u|^2\\
&+\iint_{\R^d\times\R^d}
\bigl(\nabla\phi_R(x)-\nabla\phi_R(y)\bigr)\cdot\nabla K(x-y)
|u(x)|^2|u(y)|^2\dd x\dd y.
\end{aligned}
\]
It follows first for smooth approximations and then for $H^1$ solutions in distributions in time. In this passage, the symmetrized interaction is bounded in absolute value by $C K(x-y)|u(x)|^2|u(y)|^2$, since $\|D^2\phi_R\|_\infty\le C$. Thus the Hardy--Littlewood--Sobolev estimate and $H^1$ convergence justify the limit even in dimension three, where $\nabla K$ alone is not locally integrable. To interpret the symmetrized interaction, one may first replace $K$ by $(|x|^2+\varepsilon^2)^{-1}$; the same bound is uniform in $\varepsilon$.

As $R\to\infty$, we have $D^2\phi_R\to2I$, $\|\Delta^2\phi_R\|_\infty\to0$, and the interaction integrand tends to $-4K(x-y)|u(x)|^2|u(y)|^2$. On compact time intervals these limits are dominated by the bounded $H^1$ norm and the Hartree energy. Moreover, $0\le\phi_R(x)\le C|x-x_*|^2$, so \eqref{eq5.11} permits passage to the limit on the left. Consequently,
\[
V_{x_*}''(t)=8A(u(t))-4\Dcal(u(t))=16E_0
\]
in distributions, and hence $V_{x_*}$ is a quadratic polynomial. Its nonnegativity and \eqref{eq5.11} force
\[
V_{x_*}(t)=8E_0(T-t)^2.
\]
Passing to the limit in the integrated first identity gives
\[
V_{x_*}'(t)=4\operatorname{Im}\int_{\R^d}(x-x_*)\cdot\overline u\nabla u.
\]
The right side is continuous on $[0,T)$: $u$ is continuous in $H^1$, and \eqref{eq5.11} makes $(x-x_*)u(t)$ locally bounded and weakly continuous in $L^2$, as follows by testing first against compactly supported functions and then using density. Thus $V_{x_*}\in C^1([0,T))$, including the right derivative at zero.
Since $u(0)$ now has finite variance, $f=e^{i|x-x_*|^2/(4T)}u(0)$ belongs to $H^1$, and
\[
E(f)=E_0+\frac{V_{x_*}'(0)}{8T}+\frac{V_{x_*}(0)}{8T^2}=0,
\qquad M(f)=m_d.
\]
Lemma~\ref{Lem4.1} and Theorem~\ref{Thm1.1} identify $f$ as $e^{i\theta}\lambda^{d/2}Q_d(\lambda(\cdot-y))$, proving the stated initial-data formula. The pseudoconformal symmetry and uniqueness of the $H^1$ flow give the corresponding solution.
\end{proof}

\section{The optimal stability constant}\label{Sec6}

This section proves Theorem~\ref{Thm1.6} and its radial counterpart, Theorem~\ref{Thm6.12}. Fix $d\ge3$ and write $Q=Q_d$. We retain $m_d$ in the normalization and set $p=2d/(d-1)$. The set $\mathcal N$, the quotient $\mathcal R$ and its infimum $c_*$ are defined in \eqref{eq1.8}. By Lemma~\ref{Lem5.1},
\[
\mathcal R(u)=\frac{2m_d-\Dcal(u)}{4m_d-2\mathcal C(u)},
\qquad u\in\mathcal N\setminus\mathcal O_d.
\]
In particular, $\mathcal R(u)<1/2$ if and only if $\Dcal(u)>\mathcal C(u)$.

For real functions, the Riesz bilinear form is
\[
\mathcal B(f,g)=\iint_{\R^d\times\R^d}|x-y|^{-2}f(x)g(y)\dd x\dd y.
\]
It is continuous on $L^{p/2}\times L^{p/2}$ by the Hardy--Littlewood--Sobolev inequality. It is also positive semidefinite: for Schwartz functions, its Fourier multiplier is a positive constant times $|\xi|^{2-d}$, and approximation in $L^{p/2}$ extends positivity to this space. Thus
\[
\mathcal B(f,g)^2\le\mathcal B(f,f)\mathcal B(g,g).
\]
We will also use the continuity estimate
\[
|\Dcal(u)-\Dcal(v)|
\le C\|u-v\|_p(\|u\|_p+\|v\|_p)^3.
\]
This follows by expanding the two density factors and applying H\"older's inequality and the bilinear Hardy--Littlewood--Sobolev estimate.

\begin{lemma}\label{Lem6.1}
Let $f_n\in H^1(\R^d;\C)$ satisfy $M(f_n)\to m_d$ and $A(f_n)\to m_d$. Define, for all sufficiently large $n$,
\begin{equation}\label{eq6.1}
\alpha_n=\left(\frac{m_d}{M(f_n)}\right)^{1/2},
\qquad
\rho_n=\left(\frac{M(f_n)}{A(f_n)}\right)^{1/2},
\qquad
u_n(x)=\alpha_n\rho_n^{d/2}f_n(\rho_nx).
\end{equation}
Then $u_n\in\mathcal N$, $\alpha_n,\rho_n\to1$, and
\[
\Dcal(u_n)-\Dcal(f_n)\to0,
\qquad
\mathcal C(u_n)-\mathcal C(f_n)\to0.
\]
Radial symmetry is preserved by this normalization.
\end{lemma}

\begin{proof}
The identities $M(f_\rho)=M(f)$, $A(f_\rho)=\rho^2A(f)$ and $\Dcal(f_\rho)=\rho^2\Dcal(f)$, where $f_\rho(x)=\rho^{d/2}f(\rho x)$, prove the normalization and the assertion about $\Dcal$. For the correlation, a change of variables in its integral representation gives
\[
|\mathcal C(f_\rho)-\mathcal C(f)|
\le\|f\|_p
\bigl\|\rho^{-d/2}W(\cdot/\rho)-W\bigr\|_{p'}.
\]
For $1/2\le\rho\le2$, exponential decay bounds the dilated kernels by $Ce^{-c|x|/2}$. Dominated convergence therefore makes the last norm tend to zero as $\rho\to1$. The norms $\|f_n\|_p$ are bounded, and the additional factor $\alpha_n\to1$ changes the correlation by $o(1)$, by H\"older's inequality. This proves the second limit.
\end{proof}

We next identify the local value of the quotient. For real perturbations, write
\[
\begin{aligned}
q_-[b]&=\|b\|_{H^1}^2-\mathcal B(Q^2,b^2),\\
q_+[a]&=\|a\|_{H^1}^2-\mathcal B(Q^2,a^2)-2\mathcal B(Qa,Qa).
\end{aligned}
\]
The tangent spaces for the two constraints and the exact $H^1$ distance are
\begin{equation}\label{eq6.2}
\begin{aligned}
T_+=\{a\in H^1(\R^d;\R):\ &\ip{Q}{a}=0,
\quad\ip{\nabla Q}{\nabla a}=0,\\
&\ip{a}{\partial_{x_j}Q}_{H^1}=0\quad(1\le j\le d)\},\\
T_-={}&\{b\in H^1(\R^d;\R):\ip{b}{Q}_{H^1}=0\}.
\end{aligned}
\end{equation}
The first two pairings in the definition of $T_+$ are the real $L^2$ pairings. Set
\begin{equation}\label{eq6.3}
\mu_+=\inf_{a\in T_+\setminus\{0\}}
\frac{q_+[a]}{\|a\|_{H^1}^2},
\qquad
\mu_-=\inf_{b\in T_-\setminus\{0\}}
\frac{q_-[b]}{\|b\|_{H^1}^2}.
\end{equation}

\begin{lemma}\label{Lem6.2}
Both infima in \eqref{eq6.3} are attained, and $0<\mu_\pm<1$.
\end{lemma}

\begin{proof}
Positivity follows from Lemma~\ref{Lem5.2}. The potential $V=|x|^{-2}*Q^2$ is positive, bounded and tends to zero at infinity. Multiplication by $V$ is consequently a compact quadratic form on $H^1$: local Rellich compactness treats a fixed ball, and $\sup_{|x|>R}V(x)\to0$ controls its complement. The nonlocal form is compact as well. If $a_n\rightharpoonup a$ in $H^1$, local compactness in $L^p$ and the estimate
\[
\|Q(a_n-a)\|_{L^{p/2}(|x|>R)}
\le\|Q\|_{L^p(|x|>R)}\|a_n-a\|_p
\]
show that $Qa_n\to Qa$ in $L^{p/2}$. Continuity of $\mathcal B$ then proves convergence of $\mathcal B(Qa_n,Qa_n)$.

Each tangent space is infinite dimensional. On a nonzero vector in either space the subtracted potential term is strictly positive; positivity of $\mathcal B(Qa,Qa)$ therefore gives a Rayleigh quotient strictly below one. For either sign, choose a minimizing sequence $h_n$ of unit $H^1$ norm, with weak limit $h$ in the same tangent space. Write the form as $q[h]=\|h\|_{H^1}^2-K[h]$, with $K$ compact. Then $K[h]=1-\mu>0$, and
\[
\mu\|h\|_{H^1}^2\le q[h]
=\|h\|_{H^1}^2-(1-\mu).
\]
Since $\mu<1$, this forces $\|h\|_{H^1}\ge1$. Weak lower semicontinuity gives the reverse inequality. Thus $h_n\to h$ strongly and $q[h]=\mu$.
\end{proof}

\begin{proposition}\label{Prop6.3}
The exact local value of the stability quotient is
\begin{equation}\label{eq6.4}
\lim_{\delta\downarrow0}
\inf_{\substack{u\in\mathcal N\\
0<\dist_{H^1}(u,\mathcal O_d)<\delta}}
\mathcal R(u)=2\min\{\mu_+,\mu_-\}.
\end{equation}
\end{proposition}

\begin{proof}
Choose a nearest ground state using Lemma~\ref{Lem5.1}, and apply its inverse phase and translation to $u$. Write $u=Q+h$, $h=a+ib$. Then
\[
\|h\|_{H^1}=\dist_{H^1}(u,\mathcal O_d),
\qquad
\ip{b}{Q}_{H^1}=0,
\qquad
\ip{a}{\partial_{x_j}Q}_{H^1}=0.
\]
The exact constraints give
\[
\ip{Q}{a}=-\frac12M(h),
\qquad
\ip{\nabla Q}{\nabla a}=-\frac12A(h).
\]
Set
\[
c=-\frac{M(h)}{2m_d},
\qquad
e=\frac{M(h)-A(h)}{2m_d},
\qquad
\widehat a=a-cQ-e\Lambda Q.
\]
The identities $\ip{Q}{\Lambda Q}=0$ and $\ip{\nabla Q}{\nabla\Lambda Q}=m_d$ imply $\widehat a\in T_+$. Radial symmetry preserves all the $H^1$ translation conditions, and
\[
\|\widehat a-a\|_{H^1}\le C\|h\|_{H^1}^2.
\]
The forms are bounded on $H^1$, so replacing $a$ by $\widehat a$ changes its quadratic form and squared norm by $O(\|h\|_{H^1}^3)$. Taylor expansion of the quartic action, as in the proof of Theorem~\ref{Thm1.4}, gives
\[
2m_d-\Dcal(u)
=2\bigl(q_+[a]+q_-[b]\bigr)+O(\|h\|_{H^1}^3).
\]
The remainder is uniform near $Q$, by the bilinear Hardy--Littlewood--Sobolev estimate. It follows that
\[
\mathcal R(u)\ge2\min\{\mu_+,\mu_-\}-C\|h\|_{H^1}.
\]
This proves the lower bound in \eqref{eq6.4}.

For the reverse bound, use the minimizers from Lemma~\ref{Lem6.2}. The exactly normalized imaginary and real curves computed in Propositions~\ref{Prop6.4} and \ref{Prop6.5} converge to $Q$ with quotient tending to $2\mu_-$ and $2\mu_+$, respectively. Their squared distances have positive leading coefficient, so they lie outside the orbit for small nonzero parameters. These two curves prove the upper bound.
\end{proof}

We now prove the strict upper bounds. In the fourth-order expansions below, set
\[
\tau_d=\frac{\|\Lambda Q\|_{H^1}^2}{m_d}.
\]

\begin{proposition}\label{Prop6.4}
For every nonzero real $b\in T_-$,
\[
c_*<\frac{2q_-[b]}{\|b\|_{H^1}^2}.
\]
In particular, $c_*<2\mu_-$.
\end{proposition}

\begin{proof}
Put
\[
\mu=\frac{q_-[b]}{\|b\|_{H^1}^2},
\qquad r=\frac{M(b)}{m_d},
\qquad s=\frac{A(b)}{m_d}.
\]
Lemma~\ref{Lem5.2} and $V>0$ imply $0<\mu<1$, while $r,s>0$. With $f_\lambda(x)=\lambda^{d/2}f(\lambda x)$, define
\begin{equation}\label{eq6.5}
\begin{gathered}
\alpha_\varepsilon=(1+r\varepsilon^2)^{-1/2},
\qquad
\lambda_\varepsilon=
\left(\frac{1+r\varepsilon^2}{1+s\varepsilon^2}\right)^{1/2},\\
u_\varepsilon=
\alpha_\varepsilon\bigl(Q_{\lambda_\varepsilon}
+i\varepsilon b_{\lambda_\varepsilon}\bigr).
\end{gathered}
\end{equation}
Then $u_\varepsilon\in\mathcal N$. Since $\mathcal B(Q^2,b^2)=(1-\mu)m_d(r+s)$, the interaction is given exactly by
\begin{equation}\label{eq6.6}
\frac{\Dcal(u_\varepsilon)}{2m_d}
=
\frac{1+(1-\mu)(r+s)\varepsilon^2+
\dfrac{\Dcal(b)}{2m_d}\varepsilon^4}
{(1+r\varepsilon^2)(1+s\varepsilon^2)}.
\end{equation}

To compute the distance, put $F(\lambda)=\ip{Q_\lambda}{Q}_{H^1}$. Differentiating the mass and gradient scaling identities yields
\[
F(1)=2m_d,
\qquad F'(1)=m_d,
\qquad F''(1)=m_d-\|\Lambda Q\|_{H^1}^2=m_d(1-\tau_d).
\]
The absolute correlation of $Q$ with its translates has a unique maximum at zero, by the equality case in Cauchy--Schwarz. Its Hessian there is
\[
-\bigl(\ip{\partial_{x_j}Q}{\partial_{x_k}Q}_{H^1}\bigr)_{j,k},
\]
which is negative definite. The correlation and its first two translation derivatives depend continuously on $u$ in $H^1$, uniformly in the translation. The strict maximum, the decay of the correlation at infinity and this Hessian therefore imply that, for $u$ near $Q$, every maximizing translation lies in a fixed small ball and is unique there.

For the curve \eqref{eq6.5}, write
\[
F_\varepsilon(y)=\ip{Q_{\lambda_\varepsilon}}{Q(\cdot-y)}_{H^1},
\qquad
G_\varepsilon(y)=\ip{b_{\lambda_\varepsilon}}{Q(\cdot-y)}_{H^1}.
\]
At $y=0$, radial symmetry gives $\nabla F_\varepsilon(0)=0$, while
\[
G_\varepsilon(0)=O(\varepsilon^2),
\qquad
\nabla G_\varepsilon(0)=O(1).
\]
For the first estimate, transfer the dilation to $W$ in the integral pairing and use $\ip{b}{Q}_{H^1}=0$. Differentiation of the transferred kernel is justified by its bound $C(1+|x|)^2e^{-c|x|/2}$ for $\lambda$ near one and by Cauchy--Schwarz against $b$. Thus the gradient at zero of the absolute correlation is
\[
\alpha_\varepsilon
\frac{\varepsilon^2G_\varepsilon(0)\nabla G_\varepsilon(0)}
{\bigl(F_\varepsilon(0)^2+
\varepsilon^2G_\varepsilon(0)^2\bigr)^{1/2}}
=O(\varepsilon^4).
\]
The maximizing translation is consequently $O(\varepsilon^4)$, and its contribution to the maximum is $O(\varepsilon^8)$. At zero translation, the imaginary part changes the absolute value by $O(\varepsilon^6)$. Hence
\[
\mathcal C(u_\varepsilon)
=\alpha_\varepsilon F(\lambda_\varepsilon)+O(\varepsilon^6).
\]
Expanding \eqref{eq5.1}, we obtain
\begin{equation}\label{eq6.7}
\begin{split}
\frac{\dist_{H^1}(u_\varepsilon,\mathcal O_d)^2}{m_d}
={}&(r+s)\varepsilon^2\\
&+\left[\left(\frac{\tau_d}{4}-1\right)(r-s)^2
-\frac32rs\right]\varepsilon^4+o(\varepsilon^4).
\end{split}
\end{equation}
Combining \eqref{eq6.6} and \eqref{eq6.7} gives
\begin{equation}\label{eq6.8}
\begin{gathered}
\mathcal R(u_\varepsilon)=2\mu+\beta_b\varepsilon^2+o(\varepsilon^2),\\
\beta_b=
\frac{(2-5\mu)rs-\Dcal(b)/m_d-
\dfrac{\mu \tau_d}{2}(r-s)^2}{r+s}.
\end{gathered}
\end{equation}

Positivity of the Riesz form implies
\[
\mathcal B(Q^2,b^2)^2\le\Dcal(Q)\Dcal(b)=2m_d\Dcal(b),
\qquad
\frac{\Dcal(b)}{m_d}\ge\frac{(1-\mu)^2}{2}(r+s)^2.
\]
Substitution in \eqref{eq6.8} yields
\[
\beta_b\le
-\frac{\mu(1+2\mu)rs}{r+s}
-\frac{(1-\mu)^2+\mu \tau_d}{2(r+s)}(r-s)^2<0.
\]
By \eqref{eq6.7}, the trial function is outside $\mathcal O_d$ for small nonzero $\varepsilon$. Its quotient is strictly below $2\mu$. Taking a minimizer for $\mu_-$ proves the last assertion.
\end{proof}

\begin{proposition}\label{Prop6.5}
Let $a\in T_+$ realize $\mu_+$, and write
\[
\mu=\mu_+,
\qquad r=\frac{M(a)}{m_d},
\qquad s=\frac{A(a)}{m_d},
\qquad
J(a)=\left.\frac{\mathrm d}{\mathrm d\lambda}
\ip{a_\lambda}{Q}_{H^1}\right|_{\lambda=1}.
\]
Define
\begin{equation}\label{eq6.9}
\Gamma_a=
\frac{-4\mathcal B(Qa,a^2)+2\mu(r-s)J(a)}{m_d(r+s)}.
\end{equation}
If $\Gamma_a\ne0$, then $c_*<2\mu_+$. The same conclusion holds if $\Gamma_a=0$ and
\begin{equation}\label{eq6.10}
(2-5\mu)rs-\frac{\Dcal(a)}{m_d}
-\frac{\mu \tau_d}{2}(r-s)^2<0.
\end{equation}
\end{proposition}

\begin{proof}
The derivative defining $J(a)$ exists for every $a\in H^1$: transferring the dilation to $W$ gives
\[
J(a)=-\int_{\R^d}a(x)
\left(\frac d2W(x)+x\cdot\nabla W(x)\right)\dd x.
\]
The exponentially decaying differentiated kernels, multiplied by a fixed polynomial, belong to $L^2$ and justify the derivatives by Cauchy--Schwarz.

Set $u_\varepsilon=\alpha_\varepsilon(Q_{\lambda_\varepsilon}
+\varepsilon a_{\lambda_\varepsilon})$, with the parameters in \eqref{eq6.5} and the present $r,s$. The two $L^2$ tangent constraints make $M(Q+\varepsilon a)$ and $A(Q+\varepsilon a)$ exactly quadratic, so $u_\varepsilon\in\mathcal N$. Expanding the interaction and using $\ip{a}{Q}_{H^1}=0$ gives
\[
\begin{split}
\Dcal(Q+\varepsilon a)
={}&2m_d+2(1-\mu)m_d(r+s)\varepsilon^2\\
&+4\mathcal B(Qa,a^2)\varepsilon^3+\Dcal(a)\varepsilon^4.
\end{split}
\]
After normalization this expression is divided by $(1+r\varepsilon^2)(1+s\varepsilon^2)$.

The translation derivative of the correlation at zero is $O(\varepsilon^3)$: the term from $Q_{\lambda_\varepsilon}$ is zero, and the term from $a$ is zero at $\lambda=1$ by \eqref{eq6.2}, with error $O(|\lambda_\varepsilon-1|)$. The maximum argument in Proposition~\ref{Prop6.4} shows that optimizing in translation changes the correlation by $O(\varepsilon^6)$. The maximizing correlation is positive for this real curve. Since
\[
\ip{a_{\lambda_\varepsilon}}{Q}_{H^1}
=\frac{r-s}{2}J(a)\varepsilon^2+O(\varepsilon^4),
\]
the distance expansion is
\[
\begin{split}
\dist_{H^1}(u_\varepsilon,\mathcal O_d)^2
={}&m_d(r+s)\varepsilon^2-(r-s)J(a)\varepsilon^3\\
&+m_d\left[\left(\frac{\tau_d}{4}-1\right)(r-s)^2
-\frac32rs\right]\varepsilon^4+o(\varepsilon^4).
\end{split}
\]
Dividing the deficit expansion by this expression and using \eqref{eq6.9}, we obtain
\[
\mathcal R(u_\varepsilon)=2\mu+\Gamma_a\varepsilon+O(\varepsilon^2).
\]
If $\Gamma_a\ne0$, one choice of the sign of $\varepsilon$ makes the quotient strictly smaller than $2\mu$. If $\Gamma_a=0$, the coefficient of $\varepsilon^2$ is
\[
\frac{(2-5\mu)rs-\Dcal(a)/m_d-
\dfrac{\mu \tau_d}{2}(r-s)^2}{r+s},
\]
which is negative under \eqref{eq6.10}. The positive quadratic distance coefficient ensures that these trial functions are outside the orbit.
\end{proof}

\begin{lemma}\label{Lem6.6}
If $u_n\in\mathcal N$ satisfies
\[
\sup_{y\in\R^d}\int_{B_1(y)}|u_n|^2\longrightarrow0,
\]
then $\mathcal R(u_n)\to1/2$.
\end{lemma}

\begin{proof}
The vanishing lemma \cite{Lions1984} gives $u_n\to0$ in $L^p$, since $2<p<2d/(d-2)$ and the $H^1$ norms are bounded. Thus $\Dcal(u_n)\to0$ by \eqref{eq1.3}, and $\mathcal C(u_n)\to0$ by \eqref{eq5.2}. Formula \eqref{eq5.1} proves the assertion.
\end{proof}

\begin{proposition}\label{Prop6.7}
For every $d\ge3$,
\[
0<c_*\le c_*^{\mathrm{rad}}<\frac12.
\]
\end{proposition}

\begin{proof}
Positivity follows from Theorem~\ref{Thm1.4}. Fix $\lambda>1$ and put $v_\lambda=\lambda^{-1}Q_\lambda$. Then
\[
M(v_\lambda)=m_d\lambda^{-2},
\qquad A(v_\lambda)=m_d,
\qquad \Dcal(v_\lambda)=2m_d\lambda^{-2}.
\]
Choose a nonnegative radial $g\in C_c^\infty(\{1<|x|<2\})$ with $M(g)=1$, and set
\[
g_n(x)=\bigl(m_d(1-\lambda^{-2})\bigr)^{1/2}
n^{-d/2}g(x/n),
\qquad f_n=v_\lambda+g_n.
\]
Here $M(g_n)=m_d(1-\lambda^{-2})$, $A(g_n)=O(n^{-2})$ and $\|g_n\|_p=O(n^{-1/2})$. The mass cross term tends to zero because $g_n$ is supported outside $B_n$ and the $L^2$ tail of $v_\lambda$ tends to zero. The gradient cross term is $O(n^{-1})$ by Cauchy--Schwarz. Consequently $M(f_n),A(f_n)\to m_d$. Continuity in $L^p$ of the interaction and correlation gives
\[
\Dcal(f_n)\to2m_d\lambda^{-2},
\qquad
\mathcal C(f_n)\to\lambda^{-1}\mathcal C(Q_\lambda).
\]
Normalize $f_n$ by Lemma~\ref{Lem6.1}. The resulting $u_n$ are radial, and the two limits are unchanged.

Both $Q_\lambda$ and $W$ are nonnegative radial decreasing functions. For $W$, this follows because convolution preserves radial monotonicity of the positive decreasing kernel and $Q^2$, and $W=VQ$. The layer-cake representation, or the rearrangement inequality, therefore places the maximum of their convolution at zero. Hence
\begin{equation}\label{eq6.11}
\mathcal R(u_n)\longrightarrow
\Phi(\lambda)=
\frac{2m_d(1-\lambda^{-2})}
{4m_d-2\lambda^{-1}\ip{Q_\lambda}{Q}_{H^1}}.
\end{equation}
The denominator is positive, since
\[
\mathcal C(v_\lambda)
\le\sqrt{2m_d}\,\|v_\lambda\|_{H^1}
=m_d\sqrt{2(1+\lambda^{-2})}<2m_d.
\]
Finally,
\[
0<\ip{Q_\lambda}{Q}_{H^1}
=\lambda^{-d/2}\int_{\R^d}Q(z)W(z/\lambda)\dd z
\le\lambda^{-d/2}\|Q\|_1\|W\|_\infty.
\]
Thus
\[
\Phi(\lambda)=\frac12-\frac1{2\lambda^2}
+O(\lambda^{-1-d/2})\qquad(\lambda\to\infty).
\]
Since $1+d/2>2$, choose a fixed large $\lambda$ for which $\Phi(\lambda)<1/2$, and then let $n\to\infty$ in \eqref{eq6.11}. This proves $c_*^{\mathrm{rad}}\le\Phi(\lambda)<1/2$ and hence the assertion.
\end{proof}

By Lemma~\ref{Lem6.6}, minimizing sequences for either constant cannot vanish.
We now turn to minimizing sequences that stay away from the orbit. An admissible family consists of at most countably many nonzero functions $v^j\in H^1(\R^d;\C)$ such that
\begin{equation}\label{eq6.12}
\sum_jM(v^j)\le m_d,
\qquad
\sum_jA(v^j)\le m_d,
\qquad
\sup_j\mathcal C(v^j)<2m_d.
\end{equation}
The empty family is allowed, with all three quantities equal to zero. Define
\begin{equation}\label{eq6.13}
\mathcal R_{\mathrm{rel}}((v^j))
=\frac{2m_d-\sum_j\Dcal(v^j)}
{4m_d-2\sup_j\mathcal C(v^j)}.
\end{equation}
The numerator is strictly positive. Indeed, \eqref{eq1.5} implies
\[
\sum_j\Dcal(v^j)
\le\frac2{m_d}\sum_jM(v^j)A(v^j)
\le\frac2{m_d}\left(\sum_jM(v^j)\right)
\left(\sum_jA(v^j)\right)\le2m_d.
\]
Equality throughout would require a single nonzero profile with both constraints saturated and equality in the sharp inequality. Theorem~\ref{Thm1.1} and Lemma~\ref{Lem4.1} would then put it in $\mathcal O_d$, contrary to \eqref{eq6.12}. The series of interactions is therefore convergent, and the relaxed quotient is well defined and positive.

\begin{proposition}\label{Prop6.8}
Let $u_n\in\mathcal N$. After passage to a subsequence, there are at most countably many nonzero profiles $v^j\in H^1$ and translations $y_n^j$ with $|y_n^j-y_n^k|\to\infty$ for $j\ne k$, such that, for every fixed $J$,
\begin{equation}\label{eq6.14}
\begin{gathered}
u_n=\sum_{j=1}^Jv^j(\cdot-y_n^j)+r_n^J,
\qquad
\lim_{J\to\infty}\limsup_{n\to\infty}\|r_n^J\|_p=0,\\
M(u_n)=\sum_{j=1}^JM(v^j)+M(r_n^J)+o_n(1),\\
A(u_n)=\sum_{j=1}^JA(v^j)+A(r_n^J)+o_n(1).
\end{gathered}
\end{equation}
In particular, the two sum constraints in \eqref{eq6.12} hold, and
\begin{equation}\label{eq6.15}
\Dcal(u_n)\longrightarrow\sum_j\Dcal(v^j),
\qquad
\mathcal C(u_n)\longrightarrow\sup_j\mathcal C(v^j).
\end{equation}
If the distance to $\mathcal O_d$ is bounded away from zero, this family is admissible and
\[
\mathcal R(u_n)\longrightarrow\mathcal R_{\mathrm{rel}}((v^j)).
\]
\end{proposition}

\begin{proof}
Use the translation profile decomposition for bounded $H^1$ sequences from \cite{HmidiKeraani2005}. It applies to real and imaginary parts simultaneously and gives \eqref{eq6.14}, together with $r_n^J(\cdot+y_n^j)\rightharpoonup0$ in $H^1$ for each $j\le J$. In this decomposition only translations are extracted; the remainder may still carry mass or gradient norm.

Fix $J$ and write $v_n^j=v^j(\cdot-y_n^j)$. For each $j\le J$, the translated remainder $r_n^J(\cdot+y_n^j)$ converges to zero locally in $L^p$ and almost everywhere. Its $L^p$ norm is bounded. H\"older's inequality on a ball and the tail bound
\[
\|v^j r_n^J(\cdot+y_n^j)\|_{L^{p/2}(|x|>R)}
\le\|v^j\|_{L^p(|x|>R)}\|r_n^J\|_p
\]
therefore imply $v_n^j\overline{r_n^J}\to0$ in $L^{p/2}$. Moreover, $|r_n^J(\cdot+y_n^j)|^2\rightharpoonup0$ in $L^{p/2}$. The Hardy--Littlewood--Sobolev inequality places $|x|^{-2}*|v^j|^2$ in $L^d$, the dual space of $L^{p/2}$. Testing this weak convergence against that fixed potential shows that
\[
\mathcal B(|v_n^j|^2,|r_n^J|^2)\longrightarrow0.
\]

For two distinct profiles, put $f=v^j$, $g=v^k$ and $z_n=y_n^k-y_n^j$, so $|z_n|\to\infty$. For $R>0$, let $f_R=f\mathbf1_{B_R}$ and $g_R=g\mathbf1_{B_R}$. By the bilinear Hardy--Littlewood--Sobolev estimate and translation invariance,
\[
\begin{aligned}
&\bigl|\mathcal B(|f|^2,|g(\cdot-z_n)|^2)
-\mathcal B(|f_R|^2,|g_R(\cdot-z_n)|^2)\bigr|\\
&\qquad\le C\bigl(\|f-f_R\|_p^2\|g\|_p^2
+\|f\|_p^2\|g-g_R\|_p^2\bigr).
\end{aligned}
\]
The right side tends to zero as $R\to\infty$, uniformly in $n$. When $|z_n|>2R$, the supports of the truncated functions are separated by at least $|z_n|-2R$, and hence
\[
0\le\mathcal B(|f_R|^2,|g_R(\cdot-z_n)|^2)
\le\frac{M(f)M(g)}{(|z_n|-2R)^2}.
\]
Letting first $n\to\infty$ for fixed $R$, then $R\to\infty$, proves $\mathcal B(|v_n^j|^2,|v_n^k|^2)\to0$. The same truncations have disjoint supports for large $n$, so
\[
\|f\overline{g(\cdot-z_n)}\|_{p/2}
\le\|f-f_R\|_p\|g\|_p+\|f\|_p\|g-g_R\|_p.
\]
The same order of limits gives $v_n^j\overline{v_n^k}\to0$ in $L^{p/2}$.

It follows that, for fixed $J$,
\[
\left\||u_n|^2-\sum_{j=1}^J|v_n^j|^2-|r_n^J|^2\right\|_{p/2}
\longrightarrow0.
\]
Expanding $\mathcal B(|u_n|^2,|u_n|^2)$, the terms containing this difference vanish by the bilinear estimate and the bounded $L^{p/2}$ norms of the densities. The density cross terms vanish by the estimates above. Thus
\[
\Dcal(u_n)=\sum_{j=1}^J\Dcal(v^j)+\Dcal(r_n^J)+o_n(1).
\]
The two norm constraints and \eqref{eq1.5} imply $\sum_j\Dcal(v^j)\le2m_d$. Since $\Dcal(r_n^J)\le C\|r_n^J\|_p^4$, we obtain
\[
\limsup_{n\to\infty}\left|\Dcal(u_n)-\sum_{j=1}^\infty\Dcal(v^j)\right|
\le C\limsup_{n\to\infty}\|r_n^J\|_p^4
+\sum_{j>J}\Dcal(v^j).
\]
Both terms on the right tend to zero as $J\to\infty$, by \eqref{eq6.14} and convergence of the series. This proves the first limit in \eqref{eq6.15}.

For the second, put $F_j(y)=\ip{v^j}{Q(\cdot-y)}_{H^1}$. By Lemma~\ref{Lem5.1}, these are continuous functions vanishing at infinity. For a fixed finite family,
\[
\sup_y\left|\sum_{j=1}^JF_j(y-y_n^j)\right|
\longrightarrow\max_{1\le j\le J}\|F_j\|_\infty.
\]
To prove the upper bound, fix $\eta>0$ and choose a radius outside which every $|F_j|$ is smaller than $\eta/J$. For large $n$ the translated balls are disjoint, so at any point at most one summand is not covered by this bound. For the lower bound, evaluate at a maximum point of a profile with largest correlation; all other summands tend to zero. By \eqref{eq5.2}, the remainder changes this supremum by at most $\|W\|_{p'}\|r_n^J\|_p$. Moreover, $\mathcal C(v^j)\le\sqrt{2m_d}\|v^j\|_{H^1}\to0$ as $j\to\infty$, by norm decoupling. First let $n\to\infty$ and then $J\to\infty$. This proves the correlation limit.

If the distances are at least $\delta>0$, \eqref{eq5.1} implies $\sup_j\mathcal C(v^j)\le2m_d-\delta^2/2$. The denominator remains positive, and the last assertion follows from \eqref{eq6.15}.
\end{proof}

\begin{lemma}\label{Lem6.9}
Every admissible family is realized as the limit of the quotients of a sequence in $\mathcal N\setminus\mathcal O_d$. Consequently,
\[
c_* =\inf\{\mathcal R_{\mathrm{rel}}((v^j)):(v^j)\text{ is admissible}\}.
\]
If the family consists of one radial profile, the realizing sequence can be chosen radial and converges weakly to that profile. The empty family also has a radial realizing sequence with weak limit zero.
\end{lemma}

\begin{proof}
First consider a finite family $v^1,\ldots,v^J$. Put
\[
M_0=m_d-\sum_{j=1}^JM(v^j),
\qquad
A_0=m_d-\sum_{j=1}^JA(v^j).
\]
Choose real compactly supported smooth functions $f,g$ with $M(f)=1$ and $A(g)=1$. Place the profiles at pairwise separating centers $y_n^j$, and add the two functions
\begin{equation}\label{eq6.16}
\begin{aligned}
p_n(x)&=M_0^{1/2}n^{-d/2}f((x-z_n)/n),\\
q_n(x)&=A_0^{1/2}n^{(d-2)/2}g(n(x-z_n')).
\end{aligned}
\end{equation}
A function with zero coefficient is omitted. These satisfy
\[
\begin{aligned}
M(p_n)&=M_0,& A(p_n)&=O(n^{-2}),&\|p_n\|_p&=O(n^{-1/2}),\\
M(q_n)&=O(n^{-2}),& A(q_n)&=A_0,&\|q_n\|_p&=O(n^{-1/2}).
\end{aligned}
\]
Choose the centers so that the two supports are disjoint and lie outside growing balls around the translated profiles. All mass and gradient cross terms then tend to zero. More explicitly, approximate the finite profiles by smooth compactly supported functions with arbitrarily small $H^1$ error, separate these supports from one another and from the added functions, and bound the errors by Cauchy--Schwarz and the uniform $H^1$ bounds.

For
\[
h_n=\sum_{j=1}^Jv^j(\cdot-y_n^j)+p_n+q_n,
\]
we obtain $M(h_n),A(h_n)\to m_d$. The finite-profile calculation in Proposition~\ref{Prop6.8} and the vanishing $L^p$ norms of $p_n,q_n$ give
\[
\Dcal(h_n)\to\sum_{j=1}^J\Dcal(v^j),
\qquad
\mathcal C(h_n)\to\max_{1\le j\le J}\mathcal C(v^j).
\]
Normalize $h_n$ by \eqref{eq6.1}. Lemma~\ref{Lem6.1} shows that the parameters tend to one and that these limits are unchanged. Admissibility makes the limiting squared distance positive, so the normalized functions lie outside the orbit for large $n$.

For a countable family, apply this construction to its first $J$ profiles and fill the two corresponding deficits. For each $J$, choose an index so large that the normalized interaction and correlation differ from their finite-family limits by at most $J^{-1}$. The interaction series converges by \eqref{eq1.5}, and the maxima of the first $J$ correlations increase to the full supremum. A diagonal sequence therefore has quotient tending to \eqref{eq6.13}. The construction with no fixed profile realizes the empty family.

For a single radial profile $v$, choose $f,g$ radial and supported in $\{1<|x|<2\}$, and set both centers in \eqref{eq6.16} equal to zero. The support of $p_n$ escapes to infinity and that of $q_n$ shrinks to the origin. In addition to their small gradient and mass norms, respectively, the bounds
\[
\left|\int v\overline{p_n}\right|
\le\|v\|_{L^2(|x|>n)}\|p_n\|_2\to0,
\qquad
\left|\int\nabla v\cdot\nabla\overline{q_n}\right|
\le\|\nabla v\|_{L^2(B_{2/n})}\|\nabla q_n\|_2\to0
\]
prove the norm limits. The two added functions converge weakly to zero in $H^1$, by the same tail and absolute-continuity estimates against a fixed test function. Exact normalization preserves this weak limit: transfer its dilation to the fixed $L^2$ test functions in the mass and gradient pairings and use strong continuity of dilation there. This proves the radial assertion, including $v=0$.

The realizing sequences imply that $c_*$ is no larger than any relaxed value. Conversely, every $u\in\mathcal N\setminus\mathcal O_d$ is an admissible singleton and has the same original and relaxed quotient. Taking infima proves equality.
\end{proof}

\begin{lemma}\label{Lem6.10}
Every minimizing admissible family for \eqref{eq6.13} saturates at least one of the two sum constraints in \eqref{eq6.12}.
\end{lemma}

\begin{proof}
For an admissible family put $D_0=\sum_j\Dcal(v^j)$ and $C_0=\sup_j\mathcal C(v^j)$. Multiplication of every profile by $t>0$ changes the quotient to
\[
f(t)=\frac{2m_d-t^4D_0}{4m_d-2tC_0}
\]
whenever the scaled family remains admissible. If $D_0=0$, this function is nondecreasing and $f(1)\ge1/2$. Such a family cannot minimize, by Proposition~\ref{Prop6.7} and Lemma~\ref{Lem6.9}.

Suppose $D_0>0$. At every critical point with positive denominator, differentiation gives
\[
f''(t)=-\frac{12t^2D_0}{4m_d-2tC_0}<0.
\]
Thus $f$ has no interior local minimum. If both sum constraints were strict at a minimizing family, an open interval of $t$ containing one would preserve them and the strict correlation bound. Minimality would make $t=1$ a local minimum, a contradiction.
\end{proof}

\begin{lemma}\label{Lem6.11}
Let $v\in H^1(\R^d;\R)$ satisfy $v\ge0$, $M(v)=m_d$ and $A(v)<m_d$, and put
\[
\kappa=\left(\frac{m_d-A(v)}{m_d}\right)^{1/2}.
\]
\begin{enumerate}
\item[(a)] For every $k\in\R^d$ with $|k|=\kappa$, the function $u=e^{ik\cdot x}v$ belongs to $\mathcal N$ and $\mathcal C(u)<\mathcal C(v)$.
\item[(b)] If $v$ is radial, then $u=e^{i\kappa|x|}v$ is radial, belongs to $\mathcal N$, and $\mathcal C(u)<\mathcal C(v)$.
\end{enumerate}
In both cases $\Dcal(u)=\Dcal(v)$, and
\[
\mathcal R(u)<\frac{2m_d-\Dcal(v)}{4m_d-2\mathcal C(v)}.
\]
\end{lemma}

\begin{proof}
The multipliers in both parts belong to $W^{1,\infty}$. With $\phi(x)=k\cdot x$ or $\phi(x)=\kappa|x|$, the weak product rule gives
\[
\nabla(e^{i\phi}v)=e^{i\phi}(\nabla v+iv\nabla\phi).
\]
The cross term in its squared modulus is zero because $v$ is real. Since $|\nabla\phi|=\kappa$ almost everywhere,
\[
M(u)=M(v)=m_d,
\qquad
A(u)=A(v)+\kappa^2M(v)=m_d.
\]
The modulus, and hence the interaction, is unchanged.

For each $y$, the weight $v(x)W(x-y)$ is positive on the positive-measure set $\{v>0\}$. A level set of $e^{ik\cdot x}$ is a countable union of parallel hyperplanes, and a level set of $e^{i\kappa|x|}$ is a countable union of spheres. Both have Lebesgue measure zero. Thus the phase is not almost everywhere constant on $\{v>0\}$, and strictness in the integral triangle inequality gives
\[
\left|\int e^{i\phi(x)}v(x)W(x-y)\dd x\right|
<\int v(x)W(x-y)\dd x\le\mathcal C(v).
\]
Evaluate this strict inequality at a point attaining the left-hand supremum, which exists by Lemma~\ref{Lem5.1}. This proves $\mathcal C(u)<\mathcal C(v)$.

Finally, \eqref{eq1.5} gives $\Dcal(v)\le2A(v)<2m_d$, while Cauchy--Schwarz gives $\mathcal C(v)\le\sqrt{2m_d}\|v\|_{H^1}<2m_d$. The numerator is positive, and decreasing the correlation strictly increases the denominator. This proves the quotient comparison.
\end{proof}

\begin{proof}[Proof of Theorem~\ref{Thm1.6}]
Part (i) is Proposition~\ref{Prop6.7}, part (ii) is Proposition~\ref{Prop6.4}, and part (iii) is Proposition~\ref{Prop6.3}.

Let $u_n$ minimize $c_*$. Pass to a subsequence on which its distance to $\mathcal O_d$ either tends to zero or is bounded below by a positive number. In the first case, \eqref{eq6.4} and the definition of $c_*$ imply
\[
c_*=2\min\{\mu_+,\mu_-\}.
\]
The strict inequality $c_*<2\mu_-$ forces $c_*=2\mu_+<2\mu_-$. In the second case, Proposition~\ref{Prop6.8} produces an admissible family whose relaxed value is $c_*$. This family is nonempty, since the empty value is $1/2$. Lemma~\ref{Lem6.9} makes it a relaxed minimizer, and Lemma~\ref{Lem6.10} proves saturation of at least one sum constraint.

If the family has one profile $v$ with $M(v)=A(v)=m_d$, the two norm identities in \eqref{eq6.14} imply $r_n^1\to0$ in $H^1$. Thus $u_n(\cdot+y_n^1)\to v$ strongly. Conversely, suppose a translated, phase-rotated subsequence converges strongly. Passing to a further subsequence fixes the phase limit. Every nonzero profile center must stay a bounded distance from the center of this strong limit; otherwise its translated weak limit is zero, by strong convergence and decay of correlations of fixed $L^2$ functions. Pairwise separation allows only one such profile, and strong convergence gives it the full mass and gradient norm. This proves the equivalence in (iv). A normalized optimizer is itself an admissible minimizing singleton, and a relaxed minimizing singleton with both equalities has the original quotient $c_*$. If $c_*$ is not attained and the second alternative occurs, the minimizing family therefore either has at least two profiles or has one profile with a strict mass or gradient deficit. This proves the nonattainment alternative in (iv).

For (v), remove the constant phase and apply Lemma~\ref{Lem6.11}(a). It produces a normalized function with quotient strictly smaller than the alleged relaxed minimum, contradicting Lemma~\ref{Lem6.9}.
\end{proof}

Finally, we restrict to radial functions. Set $T_\pm^{\mathrm{rad}}=T_\pm\cap H^1_{\mathrm{rad}}(\R^d;\R)$ and define $\mu_\pm^{\mathrm{rad}}$ by restricting \eqref{eq6.3} to these spaces. Thus the translation conditions in $T_+^{\mathrm{rad}}$ are automatic, but the mass and gradient conditions remain distinct. The radial relaxed problem is
\begin{equation}\label{eq6.17}
\inf_{\substack{v\in H^1_{\mathrm{rad}}(\R^d;\C)\\
M(v)\le m_d,\ A(v)\le m_d,\ \mathcal C(v)<2m_d}}
\frac{2m_d-\Dcal(v)}{4m_d-2\mathcal C(v)}.
\end{equation}
The correlation here still contains the supremum over all translations.

\begin{theorem}\label{Thm6.12}
Let $d\ge3$. Then
\[
c_*\le c_*^{\mathrm{rad}}
<\min\bigl\{\tfrac12,2\mu_-^{\mathrm{rad}}\bigr\}.
\]
Every radial minimizing sequence has a subsequence that either converges to $\{e^{i\theta}Q_d:\theta\in\R\}$, which is possible only if $c_*^{\mathrm{rad}}=2\mu_+^{\mathrm{rad}}<2\mu_-^{\mathrm{rad}}$, or converges weakly in $H^1$, without translations, to a nonzero radial function $v$ with $M(v)\le m_d$, $A(v)\le m_d$, at least one equality, and
\[
\frac{2m_d-\Dcal(v)}{4m_d-2\mathcal C(v)}=c_*^{\mathrm{rad}}.
\]
In the latter case, if $v$ has constant phase and $M(v)=m_d$, then $A(v)=m_d$. Moreover, the convergence is strong if and only if $M(v)=A(v)=m_d$. Consequently, if $c_*^{\mathrm{rad}}$ is not attained, then either $c_*^{\mathrm{rad}}=2\mu_+^{\mathrm{rad}}<2\mu_-^{\mathrm{rad}}$, or the radial relaxed problem \eqref{eq6.17} has a minimizer with $M(v)<m_d$ or $A(v)<m_d$.
\end{theorem}

\begin{proof}
The proof of Lemma~\ref{Lem6.2} applies to the closed radial tangent spaces, so $0<\mu_\pm^{\mathrm{rad}}<1$ and both values are attained. Restricting the test curve of Proposition~\ref{Prop6.4} to a radial minimizing imaginary direction proves $c_*^{\mathrm{rad}}<2\mu_-^{\mathrm{rad}}$. Proposition~\ref{Prop6.7} gives the strict bound by $1/2$, while $c_*\le c_*^{\mathrm{rad}}$ follows from inclusion of the admissible sets.

The radial local value is $2\min\{\mu_+^{\mathrm{rad}},\mu_-^{\mathrm{rad}}\}$. To check the nearest-point issue, let radial $u_n$ have distance tending to zero from the full orbit, and choose phases and centers realizing that distance. Reflection invariance of $u_n$ gives
\[
\|Q(\cdot-y_n)-Q(\cdot+y_n)\|_{H^1}
\le2\dist_{H^1}(u_n,\mathcal O_d).
\]
The left side can tend to zero only if $y_n\to0$: it stays positive on every compact annulus of centers and tends to $\sqrt{4m_d}$ when $|y_n|\to\infty$. Thus $u_n$ approaches the centered phase orbit. For a sufficiently close radial perturbation of $Q$, the maximizing translation is unique by the maximum argument in Proposition~\ref{Prop6.4}. Rotational invariance forces that maximum to occur at zero. The nearest-point perturbation is therefore radial. The constraint correction in Proposition~\ref{Prop6.3} preserves radial symmetry, and the real and imaginary trial curves remain radial. This proves the asserted radial local value.

A bounded radial $H^1$ sequence is compact in $L^p$, by the radial compactness theorem of Strauss \cite{Strauss1977}. In the present exponent, its tail control follows directly from the radial pointwise bound: uniformly for such a sequence,
\[
\int_{|x|>R}|u_n|^p
\le C R^{-(d-1)(p-2)/2}\|u_n\|_2^2
\le C R^{-1}.
\]
Local Rellich compactness and this tail estimate prove strong $L^p$ convergence along any weakly convergent subsequence. Hence both $\Dcal$ and $\mathcal C$ are continuous under radial weak $H^1$ convergence, by their $L^p$ estimates.

By the radial recovery assertion of Lemma~\ref{Lem6.9}, every competitor in \eqref{eq6.17} is the limit of quotients of radial normalized functions. Conversely, each radial normalized function outside the orbit is a competitor for \eqref{eq6.17}. The value of this relaxed problem is therefore $c_*^{\mathrm{rad}}$.

Now let $u_n$ be a radial minimizing sequence and pass to a weak limit $v$. If its distance to $\mathcal O_d$ tends to zero along a subsequence, the preceding nearest-point argument proves convergence to the centered phase orbit. The radial local value and $c_*^{\mathrm{rad}}<2\mu_-^{\mathrm{rad}}$ then force
\[
c_*^{\mathrm{rad}}=2\mu_+^{\mathrm{rad}}<2\mu_-^{\mathrm{rad}}.
\]
Otherwise, pass to a subsequence whose distances are bounded below. Weak lower semicontinuity gives $M(v),A(v)\le m_d$, and continuity of the correlation gives $\mathcal C(v)<2m_d$. Continuity of both terms in the quotient yields
\[
\frac{2m_d-\Dcal(v)}{4m_d-2\mathcal C(v)}=c_*^{\mathrm{rad}}.
\]
The strict bound by $1/2$ excludes $v=0$. Applying the amplitude argument of Lemma~\ref{Lem6.10} to this radial singleton proves that at least one constraint is saturated; radial recovery ensures that every admissible amplitude comparison has value at least $c_*^{\mathrm{rad}}$. If $v$ has constant phase, full mass and a strict gradient deficit, Lemma~\ref{Lem6.11}(b) excludes it as well.

Finally, weak convergence is strong exactly when the squared norm of the weak limit is $2m_d$, equivalently when both $M(v)=m_d$ and $A(v)=m_d$. A strong limit in this second case lies outside the orbit and attains $c_*^{\mathrm{rad}}$. Conversely, an attaining function is a radial relaxed minimizer with both constraints saturated. This proves all the assertions.
\end{proof}

\section*{Acknowledgments}
The authors were supported by the National Natural Science Foundation of China (12301145, 12561020, 12261107) and Yunnan Fundamental Research Projects (202401AU070123, 202601AT070048).
The authors acknowledge the use of ChatGPT (OpenAI) for literature searches and language refinement during the preparation of this work. All mathematical statements, proofs, and references were independently checked by the authors, who take full responsibility for the content of the paper.

\medskip
{\bf Data availability:} Data sharing is not applicable to this article as no new data were created or analyzed in this study.

\medskip
{\bf Conflict of interest:} The authors declare that they have no conflict of interest.


\begin{thebibliography}{99}

\bibitem{AroraRoudenko2022}
A.~K.~Arora and S.~Roudenko,
Global behavior of solutions to the focusing generalized Hartree equation,
\emph{Michigan Math. J.} \textbf{71} (2022), 619--672.

\bibitem{AtkinsonHan2012}
K.~Atkinson and W.~Han,
\emph{Spherical Harmonics and Approximations on the Unit Sphere: An Introduction},
Lecture Notes in Mathematics 2044, Springer, Heidelberg, 2012.

\bibitem{Banica2004}
V.~Banica,
Remarks on the blow-up for the Schr\"odinger equation with critical mass on a plane domain,
\emph{Ann. Scuola Norm. Sup. Pisa Cl. Sci. (5)} \textbf{3} (2004), no.~1, 139--170.

\bibitem{BianchiEgnell1991}
G.~Bianchi and H.~Egnell,
A note on the Sobolev inequality,
\emph{J. Funct. Anal.} \textbf{100} (1991), 18--24.

\bibitem{BDNS2025}
M.~Bonforte, J.~Dolbeault, B.~Nazaret, and N.~Simonov,
Stability in Gagliardo--Nirenberg--Sobolev inequalities: flows, regularity, and the entropy method,
\emph{Mem. Amer. Math. Soc.} \textbf{308} (2025), no.~1554.

\bibitem{CF2013}
E.~A.~Carlen and A.~Figalli,
Stability for a GNS inequality and the Log-HLS inequality, with application to the critical mass Keller--Segel equation,
\emph{Duke Math. J.} \textbf{162} (2013), 579--625.

\bibitem{CGK2025}
S.~Chakraborty, M.~Ghosh, and D.~Karmakar,
Existence of Bianchi--Egnell stability extremizer for the Hardy--Sobolev inequality,
arXiv:2505.07039v2, 2025.

\bibitem{CGWY2026}
P.~Chen, Q.~Guo, J.~Wei, and Y.~Yu,
Uniqueness and nondegeneracy of ground state for the $L^2$-critical boson star equation,
arXiv:2609.23771, 2026.

\bibitem{DAS2014}
P.~d'Avenia and M.~Squassina,
Soliton dynamics for the Schr\"odinger--Newton system,
\emph{Math. Models Methods Appl. Sci.} \textbf{24} (2014), 553--572.

\bibitem{DEFFL2025}
J.~Dolbeault, M.~J.~Esteban, A.~Figalli, R.~L.~Frank, and M.~Loss,
Sharp stability for Sobolev and log-Sobolev inequalities, with optimal dimensional dependence,
\emph{Camb. J. Math.} \textbf{13} (2025), 359--430.

\bibitem{DS2026}
S.~Dovetta and E.~Serra,
Existence of optimizers for the sharp stability constant in the logarithmic Sobolev inequality,
arXiv:2606.23225v1, 2026.

\bibitem{FL2013}
R.~L.~Frank and E.~Lenzmann,
Uniqueness of non-linear ground states for fractional Laplacians in $\R$,
\emph{Acta Math.} \textbf{210} (2013), 261--318.

\bibitem{FLS2016}
R.~L.~Frank, E.~Lenzmann, and L.~Silvestre,
Uniqueness of radial solutions for the fractional Laplacian,
\emph{Comm. Pure Appl. Math.} \textbf{69} (2016), 1671--1726.

\bibitem{GaoWang2017}
Y.~Gao and Z.~Wang,
Minimal mass non-scattering solutions of the focusing $L^2$-critical Hartree equations with radial data,
\emph{Discrete Contin. Dyn. Syst.} \textbf{37} (2017), 1979--2007.

\bibitem{GSW2026}
J.~G\'omez, T.~Schmid, and Y.~Wu,
Multisoliton solutions and blow up for the $L^2$-critical Hartree equation,
\emph{Arch. Ration. Mech. Anal.} \textbf{250} (2026), Paper No.~79.

\bibitem{HmidiKeraani2005}
T.~Hmidi and S.~Keraani,
Blowup theory for the critical nonlinear Schr\"odinger equations revisited,
\emph{Int. Math. Res. Not.} (2005), no.~46, 2815--2828.

\bibitem{Koenig2023}
T.~K\"onig,
On the sharp constant in the Bianchi--Egnell stability inequality,
\emph{Bull. Lond. Math. Soc.} \textbf{55} (2023), 2070--2075.

\bibitem{Koenig2025}
T.~K\"onig,
Stability for the Sobolev inequality: existence of a minimizer,
\emph{J. Eur. Math. Soc.} (2025), published online first.

\bibitem{KLR2009}
J.~Krieger, E.~Lenzmann, and P.~Rapha\"el,
On stability of pseudo-conformal blowup for $L^2$-critical Hartree NLS,
\emph{Ann. Henri Poincar\'e} \textbf{10} (2009), 1159--1205.

\bibitem{Kwong1989}
M.~K.~Kwong,
Uniqueness of positive solutions of $\Delta u-u+u^p=0$ in $\R^n$,
\emph{Arch. Ration. Mech. Anal.} \textbf{105} (1989), 243--266.

\bibitem{Lenzmann2009}
E.~Lenzmann,
Uniqueness of ground states for pseudorelativistic Hartree equations,
\emph{Anal. PDE} \textbf{2} (2009), 1--27.

\bibitem{LiZhang2009}
D.~Li and X.~Zhang,
On the classification of minimal mass blowup solutions of the focusing mass-critical Hartree equation,
\emph{Adv. Math.} \textbf{220} (2009), 1171--1192.

\bibitem{Li2022}
Z.~Li,
A priori estimates, uniqueness and non-degeneracy of positive solutions of the Choquard equation,
arXiv:2201.00368, 2022.

\bibitem{Lieb1977}
E.~H.~Lieb,
Existence and uniqueness of the minimizing solution of Choquard's nonlinear equation,
\emph{Stud. Appl. Math.} \textbf{57} (1977), 93--105.

\bibitem{Lieb1983}
E.~H.~Lieb,
Sharp constants in the Hardy--Littlewood--Sobolev and related inequalities,
\emph{Ann. of Math. (2)} \textbf{118} (1983), 349--374.

\bibitem{Lions1980}
P.-L.~Lions,
The Choquard equation and related questions,
\emph{Nonlinear Anal.} \textbf{4} (1980), 1063--1072.

\bibitem{Lions1984}
P.-L.~Lions,
The concentration-compactness principle in the calculus of variations. The locally compact case. I,
\emph{Ann. Inst. H. Poincar\'e Anal. Non Lin\'eaire} \textbf{1} (1984), 109--145.

\bibitem{LZX2024}
H.~Luo, D.~Zhang, and Y.~Xu,
Uniqueness and nondegeneracy of ground states for $-\Delta u+u=(I_\alpha\star u^2)u$ in $\R^3$ when $\alpha$ is close to $2$,
\emph{Adv. Nonlinear Anal.} \textbf{13} (2024), Paper No.~20240048.

\bibitem{MaZhao2010}
L.~Ma and L.~Zhao,
Classification of positive solitary solutions of the nonlinear Choquard equation,
\emph{Arch. Ration. Mech. Anal.} \textbf{195} (2010), 455--467.

\bibitem{MMRZ2026}
Z.~Ma, C.~Miao, M.~Rosenzweig, and J.~Zheng,
Global well-posedness and scattering for mass-critical Hartree equation,
arXiv:2609.28454, 2026.

\bibitem{MMZ2026}
Z.~Ma, C.~Miao, and J.~Zheng,
Classification of minimal-mass blowup solutions to the mass-critical Hartree equation at the $L^2$-critical level,
arXiv:2610.15460, 2026.

\bibitem{Merle1993}
F.~Merle,
Determination of blow-up solutions with minimal mass for nonlinear Schr\"odinger equations with critical power,
\emph{Duke Math. J.} \textbf{69} (1993), 427--454.

\bibitem{MXZ2009}
C.~Miao, G.~Xu, and L.~Zhao,
Global well-posedness and scattering for the mass-critical Hartree equation with radial data,
\emph{J. Math. Pures Appl. (9)} \textbf{91} (2009), 49--79.

\bibitem{MXZ2010}
C.~Miao, G.~Xu, and L.~Zhao,
On the blow-up phenomenon for the mass-critical focusing Hartree equation in $\R^4$,
\emph{Colloq. Math.} \textbf{119} (2010), 23--50.

\bibitem{MorozVS2013}
V.~Moroz and J.~Van Schaftingen,
Groundstates of nonlinear Choquard equations: existence, qualitative properties and decay asymptotics,
\emph{J. Funct. Anal.} \textbf{265} (2013), 153--184.

\bibitem{MorozVS2017}
V.~Moroz and J.~Van Schaftingen,
A guide to the Choquard equation,
\emph{J. Fixed Point Theory Appl.} \textbf{19} (2017), 773--813.

\bibitem{SW2026}
T.~Schmid and Y.~Wu,
On soliton clusters and collision blow up for the $L^2$-critical Hartree equation,
arXiv:2606.30640, 2026.

\bibitem{Seok2019}
J.~Seok,
Limit profiles and uniqueness of ground states to the nonlinear Choquard equations,
\emph{Adv. Nonlinear Anal.} \textbf{8} (2019), 1083--1098.

\bibitem{Strauss1977}
W.~A.~Strauss,
Existence of solitary waves in higher dimensions,
\emph{Comm. Math. Phys.} \textbf{55} (1977), 149--162.

\bibitem{Weinstein1983}
M.~I.~Weinstein,
Nonlinear Schr\"odinger equations and sharp interpolation estimates,
\emph{Comm. Math. Phys.} \textbf{87} (1983), 567--576.

\bibitem{Weinstein1985}
M.~I.~Weinstein,
Modulational stability of ground states of nonlinear Schr\"odinger equations,
\emph{SIAM J. Math. Anal.} \textbf{16} (1985), 472--491.

\bibitem{Xiang2016}
C.-L.~Xiang,
Uniqueness and nondegeneracy of ground states for Choquard equations in three dimensions,
\emph{Calc. Var. Partial Differential Equations} \textbf{55} (2016), Paper No.~134.

\end{thebibliography}
\end{document}